\documentclass[numreferences]{kluwer}
\usepackage{graphicx}
\usepackage{amsmath}

\newtheorem{theorem}{\indent Theorem}

\newtheorem{remark}{Remark}
\newtheorem{defn}{Definition}

\begin{document}
\begin{article}
\begin{opening}

\title{A New Methodology for the Development of Numerical Methods for the Numerical Solution of the Schr\"odinger Equation}

\author{Z.A. \surname{Anastassi}\thanks{e-mail: zackanas@uop.gr}}

\author{D.S. \surname{Vlachos}\thanks{e-mail: dvlachos@uop.gr}}

\author{T.E. \surname{Simos}\thanks{Highly Cited Researcher, Active Member of the European Academy of Sciences and
Arts. Corresponding Member of the European Academy of Sciences Corresponding Member of European Academy of Arts, Sciences and Humanities, Please use the following address for all correspondence: Dr. T.E. Simos, 10 Konitsis Street, Amfithea - Paleon Faliron, GR-175 64 Athens, Greece, Tel: 0030 210 94 20 091, e-mail: tsimos.conf@gmail.com, tsimos@mail.ariadne-t.gr}}

\institute{Laboratory of Computational Sciences, Department of Computer Science and Technology, Faculty of Sciences
and Technology, University of Peloponnese, GR-221 00 Tripolis,
Greece}

\runningtitle{A New Methodology for the Development of Numerical Methods}
\runningauthor{T.E. Simos, Z.A. Anastassi, D.S. Vlachos}

\begin{abstract}
In the present paper we introduce a new methodology for the
construction of numerical methods for the approximate  solution of
the one-dimensional  Schr\"odinger equation. The new methodology is
based on the requirement of vanishing the phase-lag and its
derivatives. The efficiency of the new methodology is proved via
error analysis and numerical applications.
\end{abstract}


\classification{PACS}{02.60, 02.70.Bf, 95.10.Ce, 95.10.Eg, 95.75.Pq}
\end{opening}

\section{Introduction}
\label{Intro}
\par
The radial Schr\"odinger equation can be written as:
\begin{equation}
y '' (x) = [l(l+1)/x^{2} + V(x) - k^{2}] y(x). \label{ieqn:1}
\end{equation}
\ \noindent Many problems in theoretical physics and chemistry,
material sciences, quantum mechanics and quantum chemistry,
electronics etc. can be express via the above boundary value problem
(see for example \cite{ix78} - \cite{hertz}).

We give the definitions of some terms of (\ref{eqn:1}):
\begin{itemize}
\item  The function $W(x) = l(l+1)/x^{2} + V(x)$ is called {\it
the effective potential}. This satisfies $W(x) \rightarrow 0$ as
$x \rightarrow \infty $ \item The quantity $k^{2}$ is a real
number denoting {\it the energy} \item The quantity $l$ is a given
integer representing the {\it angular momentum} \item $V$ is a
given function which denotes the {\it potential}.
\end{itemize}

The boundary conditions are:
\begin{equation}
y(0) = 0 \label{ieqn:2}
\end{equation}
\noindent and a second boundary condition, for large values of
$x$, determined by physical considerations.

\par
The last years an extended study on the development of numerical
methods for the solution of the Schr\"odinger equation has been
done. The aim of this research is the development of fast and
reliable methods for the solution of the Schr\"odinger equation and
related problems (see for example \cite{simos00_r} - \cite{siex99},
\cite{berghe2004} - \cite{jnaiam3_11}).
\par
We can divide the numerical methods for the approximate solution of
the Schr\"odinger equation and related problems into two main
categories:
\begin{enumerate}
\item Methods with constant coefficients \item Methods with
coefficients depending on the frequency of the problem
\footnote{When using a functional fitting algorithm for the solution
of the radial Schr\"odinger equation, the fitted frequency is equal
to: $\sqrt{|l(l+1)/x^{2} + V(x) - k^{2}|}$}.
\end{enumerate}
\noindent The purpose of this paper is to introduce a new
methodology for the construction of numerical methods for the
approximate  solution of the one-dimensional  Schr\"odinger equation
and related problems. The new methodology is based on the
requirement of vanishing the phase-lag and its derivatives. The
efficiency of the new methodology will be studied via the error
analysis and the application of the investigated methods to the
numerical solution of the radial Schr\"odinger equation.

More specifically, we will develop a family of hybrid Numerov-type
methods of sixth algebraic order. The development of the new family
is based on the requirement of vanishing the phase-lag and its
derivatives. We will investigate the stability and the error of the
methods of the new family. Finally, we will apply both categories of
methods the new obtained method to the resonance problem. This is
one of the most difficult problems arising from the radial
Schr\"odinger equation. The paper is organized as follows. In
Section 2 we present the theory of the new methodology. In section 3
we present the development of the new family of methods. The error
analysis is presented in section 4. In section 5 we will
investigate the stability properties of the new developed methods.
In Section 6 the numerical results are presented. Finally, in
Section 7 remarks and conclusions are discussed.

\section{Phase-lag analysis of symmetric multistep methods}

For the numerical solution of the initial value problem

\begin{equation}
\label{ivp_definition}
    y'' = f(x,y)
\end{equation}

\noindent consider a multistep method with $m$ steps which can be
used over the equally spaced intervals
$\left\{x_{i}\right\}^{m}_{i=0} \in [a,b]$ and
$h=|x_{i+1}-x_{i}|$, \, $i=0(1)m-1$.

If the method is symmetric then $a_i=a_{m-i}$ and $b_i=b_{m-i}$, \, $i=0(1)\lfloor \frac{m}{2} \rfloor$.

When a symmetric $2k$-step method, that is for $i=-k(1)k$, is
applied to the scalar test equation

\begin{equation}
\label{istab_eq} y''=-\omega^2 y
\end{equation}

a difference equation of the form
\begin{eqnarray}
\label{phl_multi_de} \nonumber A_{k}(H) \, y_{n + k} + \ldots +
A_{1}(H) \, y_{n + 1} + A_{0}(H)\, y_{n} + \nonumber \\ + A_{1}(H)\,
y_{n - 1} + ... + A_{k}(H)\, y_{n - k} = 0
\end{eqnarray}

\noindent is obtained, where $H = \omega h$, $h$ is the step length
and $A_{0} (H)$, $A_{1} (H),\ldots$, $ A_{k} (H)$
are polynomials of $H$.

The characteristic equation associated with (\ref{phl_multi_de}) is
given by:

\begin{eqnarray}
\label{phl_multi_ce} A_{k} (H)\, \lambda^{k} + ... + A_{1} (H)\,
\lambda + A_{0} (H) + A_{1} (H)\, \lambda^{ - 1} + ...\\
 + A_{k} (H)\, \lambda^{ - k} = 0
\end{eqnarray}

\begin{theorem}
\emph{\cite{simoscam97a}} The symmetric $2k$-step method with
characteristic equation given by (\ref{phl_multi_ce}) has phase-lag
order $q$ and phase-lag constant $c$ given by

\begin{equation}
\textstyle \label{iphl_multi_defn} - c\, H ^{q + 2} + O(H^{q + 4}) = \\
{\frac{{2\,
A_{k}(H)\, \cos (k\, H ) + \ldots + 2\, A_{j}(H)\, \cos (j\, H ) +
\ldots + A_{0}(H)}}{{2\, k^{2}\, A_{k}(H) + \ldots + 2\, j^{2}\,
A_{j}(H) + \ldots + 2\,A_{1}(H)}}}
\end{equation}
\end{theorem}

The formula proposed from the above theorem gives us a direct method
to calculate the phase-lag of any symmetric $2k$- step method.

\section{The New Family of Numerov-Type Hybrid Methods - Construction of the New Methods}

\subsection{First Method of the Family}

\noindent We introduce the following family of methods to
integrate $y'' = f(x,y)$ :

\begin{eqnarray}
\overline{y}_{n} = y_{n} - a_{0} \, h^2 \, \Bigl ( y''_{n+1} - 2 \,
y''_{n} + y''_{n-1} \Big ) \nonumber \\
y_{n+1} + c_{1} \,y_{n} + y_{n-1} = h^2 \, \left [ b_{0} \, \left (
y''_{n+1} + y''_{n-1} \right ) + b_{1} \overline{y}''_{n} \right ]
\label{nm}
\end{eqnarray}

The application of the above method to the scalar test equation
(\ref{stab_eq}) gives the following difference equation:

\begin{equation}
\label{inm1} \nonumber A_{1}(H) \, y_{n + 1} + A_{0}(H)\, y_{n} +
A_{1}(H)\, y_{n - 1} = 0
\end{equation}

\noindent where $H = \omega h$, $h$ is the step length and $A_{0}
(H)$ and $A_{1} (H)$ are polynomials of $H$.

The characteristic equation associated with (\ref{nm1}) is given by:

\begin{eqnarray}
\label{nm2} A_{1}(H)\, \lambda + A_{0}(H) + A_{1} (H)\,\lambda^{ -
1} = 0
\end{eqnarray}

\noindent where

\begin{eqnarray*}
{A_{1}(H)} = 1 + H^{2}\,{b_{0}} + H^{4}\,{b_{1}}\,{a_{0}} \\
{A_{0}(H)} = {c_{1}} + H^{2}\,{b_{1}} - 2\,H^{4}\,{b_{1}}\,{a_{0}}
\end{eqnarray*}

By applying $k=1$ in the formula (\ref{phl_multi_defn}), we have that
the phase-lag is equal to:

\begin{eqnarray}
phl = \frac{2\, A_{1}(H)\, \cos(H) + A_{0}(H)}{2\, A_{1}(H)} \nonumber \\
= {\displaystyle \frac {1}{2}} \,{\displaystyle \frac {2\,(1 +
H^{2}\,{b_{0}} + H^{4}\,{b_{1}}\,{a_{0}})\, \mathrm{cos}(H) +
{c_{1}} + H^{2}\,{b_{1}} - 2\,H^{4}\,{b_{1}}\,{ a_{0}}}{1 +
H^{2}\,{b_{0}} + H^{4}\,{b_{1}}\,{a_{0}}}} \label{nm3}
\end{eqnarray}

We demand that the phase-lag is equal to zero and we consider that:

\begin{equation}
b_{0} = \frac{1}{12}, \, \, b_{1} = \frac{5}{6}, \, \, c_{1} = - 2
\end{equation}

\noindent Then we find out that:

\begin{equation}
{a_{0}} = {\displaystyle \frac { - 12\,\mathrm{cos}(H) -
\mathrm{cos}(H)\,H^{2} + 12 - 5\,H^{2}}{10\,\mathrm{cos}(H)\,H^{4 }
- 10\,H^{4}}} \label{inm4}
\end{equation}

For small values of $|H|$ the formulae given by (\ref{nm4}) are
subject to heavy cancellations. In this case the following Taylor
series expansions should be used:

\begin{eqnarray}
a_{0} = {\displaystyle \frac {1}{200}}  + {\displaystyle \frac
{1}{5040}} \,H^{2} + {\displaystyle \frac {1}{144000}} \,H ^{4} +
{\displaystyle \frac {1}{4435200}} \,H^{6} \nonumber \\ +
{\displaystyle \frac {691}{99066240000}} \,H^{8} +
{\displaystyle \frac {1}{4790016000}} \,H^{10} \nonumber \\
+ {\displaystyle \frac {3617}{592812380160000}} \,H^{12}
 + {\displaystyle \frac {43867}{250445794959360000}} \,H^{14} \nonumber \\
+ {\displaystyle \frac {174611}{35213055381504000000}} \, H^{16} +
\ldots \label{inm4t}
\end{eqnarray}

The behavior of the coefficients is given in the following Figure
1.

\setcounter{figure}{0}
\begin{figure}[htb!]
\begin{center}
\includegraphics[width=\textwidth]{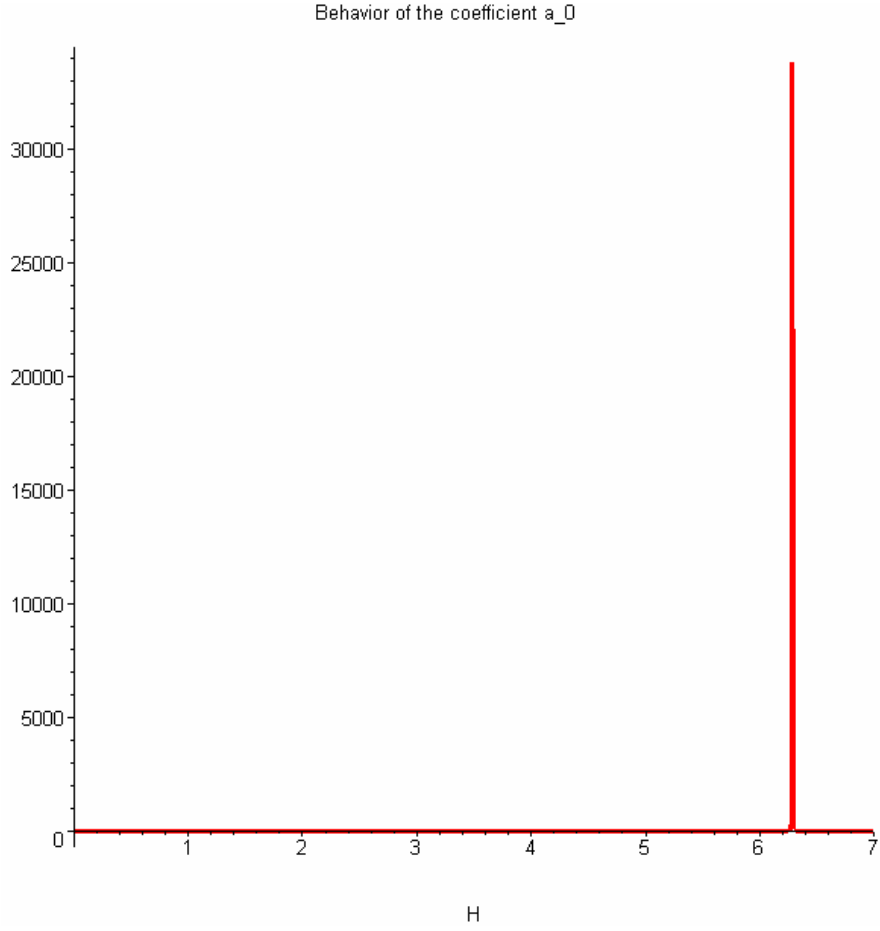}
\caption[]{\label{fig0}
Behavior of the coefficient $a_{0}$ of the new method given by
(\ref{nm4}) for several values of $H$.} \normalsize
\end{center}
\end{figure}

\vspace*{0.7cm}

The local truncation error of the new proposed method is given by:

\begin{equation}
\mathrm{LTE} = \frac{h^8}{6048} \left ( \,y_{n}^{(8)} +
\omega^{2}\,\,y_{n}^{(6)} \right ) \label{ilte}
\end{equation}

\begin{remark}
The method developed in this section is the same with the obtained
by Simos in \cite{simosijcm92}
\end{remark}

\subsection{Second Method of the Family}

\noindent Consider the family of methods presented in (\ref{nm}).

The application of the above method to the scalar test equation
(\ref{stab_eq}) gives the difference equation (\ref{nm1}) and the
characteristic equation (\ref{nm2}).

By applying $k=1$ in the formula (\ref{phl_multi_defn}), we have that
the phase-lag is given by (\ref{nm3}). The first derivative of the
phase-lag is given by:

\begin{eqnarray}
\dot{phl} = {\displaystyle \frac {1}{2}} \,{\displaystyle \frac
{T_{4} - 2\,\mathrm{T_{0}}\,\mathrm{sin}(H) + 2\,H\,{b_{1}} - 8\,H
^{3}\,{b_{1}}\,{a_{0}}}{\mathrm{T_{0}}}}  \nonumber \\
- {\displaystyle \frac {1}{2}} \,{\displaystyle \frac {(2
\,\mathrm{T_{0}}\,\mathrm{cos}(H) + {c_{1}} + H^{2}\,{b_{1}} - 2\,H
^{4}\,{b_{1}}\,{a_{0}})\,(2\,H\,{b_{0}} + 4\,H^{3}\,{b_{1}}\,{a_{
0}})}{\mathrm{T_{0}}^{2}}}  \nonumber \\
\mathrm{T_{0}} = 1 + H^{2}\,{b_{0}} + H^{4}\,{b_{1}}\,{a_{0}}
\nonumber \\
\mathrm{T_{4}} = 2\,(2\,H\,{b_{0}} +
4\,H^{3}\,{b_{1}}\,{a_{0}})\,\mathrm{ cos}(H) \label{nm3a}
\end{eqnarray}

We demand that the phase-lag and its derivative are equal to zero and we consider that:

\begin{equation}
b_{0} = \frac{1}{12}, \, \, b_{1} = \frac{5}{6}
\end{equation}

\noindent Then we find out that:

\begin{eqnarray}
{a_{0}} = {\displaystyle \frac { - \mathrm{sin}(H)\,H^{2} + 10\, H +
2\,\mathrm{cos}(H)\,H - 12\,\mathrm{sin}(H)}{10\,\mathrm{sin}
(H)\,H^{4} - 40\,\mathrm{cos}(H)\,H^{3} + 40\,H^{3}}} \nonumber
\\
{c_{1}} = (24\,\mathrm{cos}(2\,H) + 24 - 48\,\mathrm{cos}(H) + H
^{2}\,\mathrm{cos}(2\,H) \nonumber \\ - 9\,H^{2} +
8\,\mathrm{cos}(H)\,H^{2} - 6\,H^{3}\,\mathrm{sin}(H) \nonumber \\
- 12\,\mathrm{sin}(H)\,H)/(6\,\mathrm{sin}(H)\,H - 24\,
\mathrm{cos}(H) + 24)\label{nm4a}
\end{eqnarray}

For small values of $|H|$ the formulae given by (\ref{nm4a}) are
subject to heavy cancelations. In this case the following Taylor
series expansions should be used:

\begin{eqnarray}
a_{0} = {\displaystyle \frac {1}{200}}  + {\displaystyle \frac
{1}{3780}} \,H^{2} + {\displaystyle \frac {73}{5443200}} \, H^{4} +
{\displaystyle \frac {509}{769824000}} \,H^{6} \nonumber \\ +
{\displaystyle \frac {2833543}{88268019840000}} \,H^{8} \\
+ {\displaystyle \frac {4912333}{3177648714240000}} \,H^{ 10}
\nonumber \\  + {\displaystyle \frac
{288303913}{3889442026229760000}} \,H^{12} + {\displaystyle \frac
{165095552521}{ 46556621053970227200000}} \,H^{14} \nonumber \\  +
{\displaystyle \frac { 15619496804053}{92182109686861049856000000}}
\,H^{16}  + \ldots \nonumber \\
c_{1} - 2 + {\displaystyle \frac {1}{18144}} \,H^{8}
 + {\displaystyle \frac {13}{16329600}} \,H^{10} +
{\displaystyle \frac {31}{461894400}} \,H^{12} \nonumber \\   +
{\displaystyle \frac {308851}{105921623808000}} \,H^{14} +
{\displaystyle \frac {537907}{3813178457088000}} \,H^{16} + \ldots
\label{nm4ta}
\end{eqnarray}

The behavior of the coefficients is given in the following Figure 2.

\setcounter{figure}{1}
\begin{figure}[htb!]
\begin{center}
\begin{tabular}{ccc}
\includegraphics[width=0.4\textwidth]{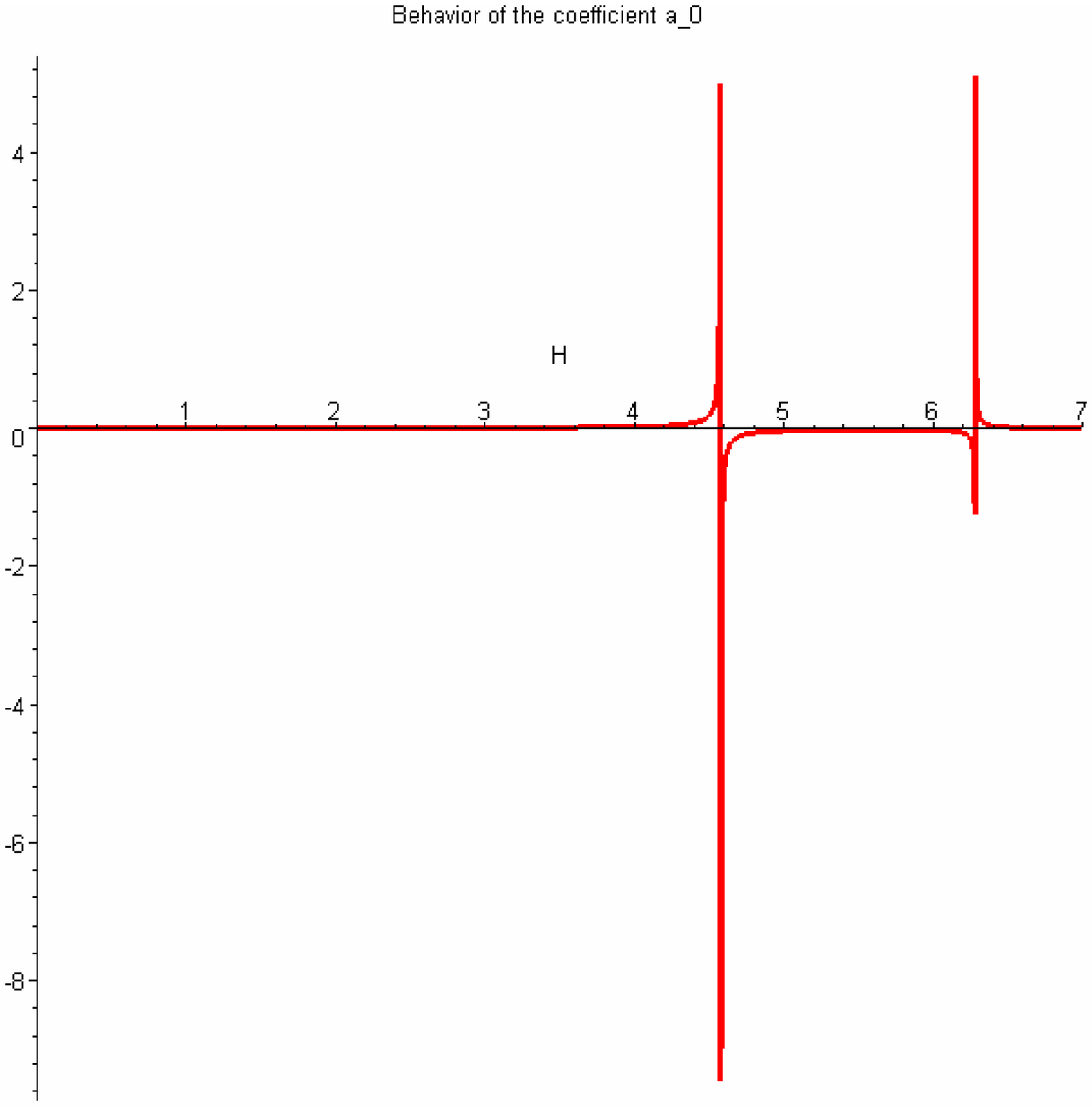}&\quad &
\includegraphics[width=0.4\textwidth]{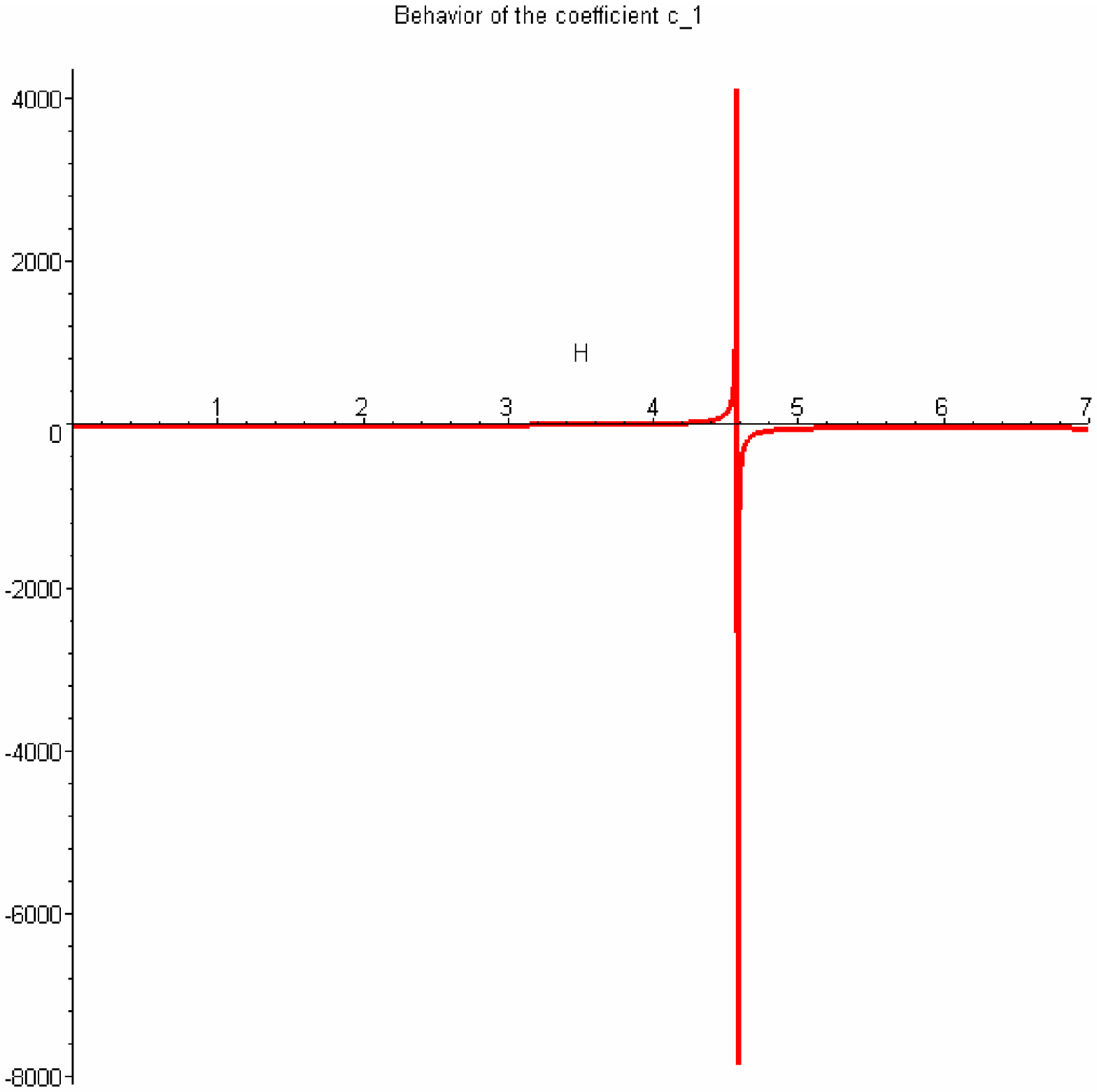}\\
\end{tabular}
\caption[]{\label{fig2}
Behavior of the coefficients of the new method given by (\ref{nm4a})
for several values of $H$.} \normalsize
\end{center}
\end{figure}

\vspace*{0.7cm}

The local truncation error of the new proposed method is given by:

\begin{equation}
\mathrm{LTE} = \frac{h^8}{18144} \left (3\, \,y_{n}^{(8)} + 4\,
\omega^{2}\,\,y_{n}^{(6)} + \omega^8 \, y_{n} \right ) \label{ltea}
\end{equation}

\subsection{Third Method of the Family}

\noindent Consider the family of methods presented in (\ref{nm}).

The application of the above method to the scalar test equation
(\ref{stab_eq}) gives the difference equation (\ref{nm1}) and the
characteristic equation (\ref{nm2}).

By applying $k=1$ in the formula (\ref{phl_multi_defn}), we have that
the phase-lag is given by (\ref{nm3}). The first derivative of the
phase-lag is given by (\ref{nm3a}). The second derivative of the
phase-lag can be written as:

\begin{eqnarray}
\ddot{phl} = {\displaystyle \frac {1}{2}} \,{\displaystyle \frac
{T_{3} - 4\,\mathrm{T_{2}}\,\mathrm{sin}(H) - 2\,\mathrm{T_{1}}\,
\mathrm{cos}(H) + 2\,{b_{1}} - 24\,{b_{1}}\,{a_{0}}\,H^{2}}{
\mathrm{T_{1}}}}  \nonumber \\
 - {\displaystyle \frac {(2\,\mathrm{T_{2}}\,\mathrm{cos}(H)
 - 2\,\mathrm{T_{1}}\,\mathrm{sin}(H) + 2\,H\,{b_{1}} - 8\,H^{3}\,{
b_{1}}\,{a_{0}})\,\mathrm{T_{2}}}{\mathrm{T_{1}}^{2}}}  \nonumber \\
 + {\displaystyle \frac {(2\,\mathrm{T_{1}}\,\mathrm{cos}(H)
 + {c_{1}} + H^{2}\,{b_{1}} - 2\,H^{4}\,{b_{1}}\,{a_{0}})\,
\mathrm{T_{2}}^{2}}{\mathrm{T_{1}}^{3}}}  \nonumber \\
 - {\displaystyle \frac {1}{2}} \,{\displaystyle \frac {(2
\,\mathrm{T_{1}}\,\mathrm{cos}(H) + {c_{1}} + H^{2}\,{b_{1}} - 2\,H
^{4}\,{b_{1}}\,{a_{0}})\,(2\,{b_{0}} + 12\,{b_{1}}\,{a_{0}}\,H^{2
})}{\mathrm{T_{1}}^{2}}}  \nonumber \\
\mathrm{T_{1}} = 1 + H^{2}\,{b_{0}} + H^{4}\,{b_{1}}\,{a_{0}} \nonumber \\
\mathrm{T_{2}} = 2\,H\,{b_{0}} + 4\,H^{3}\,{b_{1}}\,{a_{0}} \nonumber \\
\mathrm{T_{3}} = 2\,(2\,{b_{0}} +
12\,{b_{1}}\,{a_{0}}\,H^{2})\,\mathrm{cos }(H) \label{nm3b}
\end{eqnarray}

We demand that the phase-lag and its first and second derivative are equal to zero and we consider that:

\begin{equation}
b_{0} = \frac{1}{12}
\end{equation}

\noindent Then we find out that:

\begin{eqnarray}
{a_{0}} = {\displaystyle \frac {1}{2}} \Bigl (
\mathrm{cos}(H)\,H^{3}
 + 12\,\mathrm{cos}(H)\,H - 12\,\mathrm{sin}(H) + 3\,\mathrm{sin}
(H)\,H^{2} \Bigr ) / \nonumber \\  \Bigl ( \Bigl (
\mathrm{cos}(H)^{2}\,H^{3} + 16\,\mathrm{cos}(H)^{2}\,H +
5\,\mathrm{cos}(H)\,H^{2}\,\mathrm{sin}(H) \nonumber \\
+ 72\,\mathrm{cos}(H)\,\mathrm{sin}(H) + 2\,\mathrm{cos}(H)\,H^{3}
+ 32 \,\mathrm{cos}(H)\,H \nonumber \\
+ 2\,\mathrm{sin}(H)\,H^{2} - 48\,H - 2\,H^{3} -
72\,\mathrm{sin}(H) \Bigr ) H^{2} \Bigr ) \nonumber \\
{c_{1}} = {\displaystyle \frac {1}{6}} \Bigl (
24\,\mathrm{cos}(H)^{2} \,H^{2} + \mathrm{cos}(H)^{2}\,H^{4} +
96\,\mathrm{cos}(H)^{2} \nonumber \\
+ \mathrm{cos}(H)\,\mathrm{sin}(H)\,H^{3} + 12\,\mathrm{cos}(H)\,H
^{2} \nonumber \\
- 24\,\mathrm{cos}(H)\,\mathrm{sin}(H)\,H - 96\,\mathrm{ cos}(H) +
\mathrm{cos}(H)\,H^{4} \nonumber \\
- \mathrm{sin}(H)\,H^{3} - 2\,H^{4} - 60\,\mathrm{sin}(H)\,H -
48\,H^{2} \Bigr )  / \nonumber \\ \Bigl ( \mathrm{cos}(H)\,H^{2} +
7\,\mathrm{sin}(H) \,H + 8 -
8\,\mathrm{cos}(H) \Bigr ) \nonumber \\
{b_{1}} =  - {\displaystyle \frac {1}{6}} \Bigl (
\mathrm{cos}(H)^{2}\, H^{3} + 16\,\mathrm{cos}(H)^{2}\,H +
5\,\mathrm{cos}(H)\,H^{2}\,
\mathrm{sin}(H) \nonumber \\
+ 72\,\mathrm{cos}(H)\,\mathrm{sin}(H) + 2\,
\mathrm{cos}(H)\,H^{3} \nonumber \\
+ 32\,\mathrm{cos}(H)\,H + 2\,\mathrm{sin}(H)\,H^{2} - 48 \,H -
2\,H^{3} - 72\,\mathrm{sin}(H) \Bigr ) / \nonumber \\ \Bigl ( H
\Bigl ( \mathrm{cos}(H)\,H^{2} + 7\,\mathrm{sin}(H)\,H + 8 -
8\,\mathrm{ cos}(H) \Bigr ) \Bigr )\label{nm4b}
\end{eqnarray}

For small values of $|H|$ the formulae given by (\ref{nm4b}) are
subject to heavy cancellations. In this case the following Taylor
series expansions should be used:

\begin{eqnarray}
a_{0} = {\displaystyle \frac {1}{200}}  + {\displaystyle \frac
{1}{2520}} \,H^{2} + {\displaystyle \frac {31}{907200}} \,H ^{4} +
{\displaystyle \frac {1229}{1197504000}} \,H^{6} \nonumber \\ +
{\displaystyle \frac {18427}{980755776000}} \,H^{8} - {\displaystyle
\frac {669341}{98075577600000}} \,H^{10} \nonumber \\
 - {\displaystyle \frac {13764419}{25184162304000000}} \,H^{12}
 - {\displaystyle \frac {281298850211}{5747730994317312000000}} \,H^{14} \nonumber \\
 - {\displaystyle \frac {161773544323}{103459157897711616000000}} \,H^{16} +
 \ldots \nonumber \\
c_{1} = - 2 - {\displaystyle \frac {1}{6048}} \,H^{8} -
{\displaystyle \frac {17}{2721600}} \,H^{10} - {\displaystyle \frac
{43}{57480192}} \,H^{12} \nonumber \\ - {\displaystyle \frac
{1515133}{ 23538138624000}} \,H^{14} - {\displaystyle \frac
{25819}{4483454976000}} \,H^{16} + \ldots \nonumber \\
b_{1} ={\displaystyle \frac {5}{6}}  + {\displaystyle \frac
{1}{3024}} \,H^{6} + {\displaystyle \frac {11}{725760}} \,H ^{8} +
{\displaystyle \frac {2353}{1437004800}} \,H^{10} \nonumber \\ +
{\displaystyle \frac {186533}{1307674368000}} \,H^{12}  +
{\displaystyle \frac {112457}{8826801984000}} \,H^{14} \nonumber \\
 + {\displaystyle \frac {1635421}{1440534083788800}} \,H^{16} +
 \ldots \label{nm4tb}
\end{eqnarray}

The behavior of the coefficients is given in the following Figure 3.

\setcounter{figure}{2}
\vspace*{0.7cm}

The local truncation error of the new proposed method is given by:

\begin{equation}
\mathrm{LTE} = \frac{h^8}{6048} \left (y_{n}^{(8)} + 2\,
\omega^{2}\,\,y_{n}^{(6)} - 2 \omega^6 \, y_{n}^{(2)} - \omega^8 \,
y_{n} \right ) \label{lteb}
\end{equation}

\subsection{Fourth Method of the Family}

\noindent Consider the family of methods presented in (\ref{nm}).

The application of the above method to the scalar test equation
(\ref{stab_eq}) gives the difference equation (\ref{nm1}) and the
characteristic equation (\ref{nm2}).

By applying $k=1$ in the formula (\ref{phl_multi_defn}), we have that
the phase-lag is given by (\ref{nm3}). The first derivative of the
phase-lag is given by (\ref{nm3a}). The second derivative of the
phase-lag is given by (\ref{nm3b}). The third derivative of the
phase-lag can be written as:

\begin{eqnarray}
\dddot{phl} = {\displaystyle \frac {1}{2}} \,{\displaystyle \frac
{T_{9} - 6\,\mathrm{T_{8}}\,\mathrm{cos}(H) + 2\,
\mathrm{T_{5}}\,\mathrm{sin}(H) - 48\,{b_{1}}\,{a_{0}}\,H}{\mathrm{
T_{5}}}}  \nonumber \\
 - {\displaystyle \frac {3}{2}} \,{\displaystyle \frac {(2
\,\mathrm{T_{7}}\,\mathrm{cos}(H) -
4\,\mathrm{T_{8}}\,\mathrm{sin}(H ) -
2\,\mathrm{T_{5}}\,\mathrm{cos}(H) + 2\,{b_{1}} - 24\,{b_{1}}\,
{a_{0}}\,H^{2})\,\mathrm{T_{8}}}{\mathrm{T_{5}}^{2}}}  \nonumber \\
 + {\displaystyle \frac {3\,(2\,\mathrm{T_{8}}\,\mathrm{cos}
(H) - 2\,\mathrm{T_{5}}\,\mathrm{sin}(H) + 2\,H\,{b_{1}} - 8\,H^{3}
\,{b_{1}}\,{a_{0}})\,\mathrm{T_{8}}^{2}}{\mathrm{T_{5}}^{3}}}  \nonumber \\
 - {\displaystyle \frac {3}{2}} \,{\displaystyle \frac {(2
\,\mathrm{T_{8}}\,\mathrm{cos}(H) -
2\,\mathrm{T_{5}}\,\mathrm{sin}(H ) + 2\,H\,{b_{1}} -
8\,H^{3}\,{b_{1}}\,{a_{0}})\,\mathrm{T_{7}}}{ \mathrm{T_{5}}^{2}}} -
{\displaystyle \frac {3\,\mathrm{T_{6}}\,
\mathrm{T_{8}}^{3}}{\mathrm{T_{5}}^{4}}}  \nonumber \\
 + {\displaystyle \frac {3\,\mathrm{T_{6}}\,\mathrm{T_{8}}\,
\mathrm{T_{7}}}{\mathrm{T_{5}}^{3}}}  - {\displaystyle \frac {12\,
\mathrm{T_{6}}\,{b_{1}}\,{a_{0}}\,H}{\mathrm{T_{5}}^{2}}}  \nonumber \\
\mathrm{T_{5}} = 1 + H^{2}\,{b_{0}} + H^{4}\,{b_{1}}\,{a_{0}} \nonumber \\
\mathrm{T_{6}} = 2\,\mathrm{T_{5}}\,\mathrm{cos}(H) + {c_{1}} + H^{2
}\,{b_{1}} - 2\,H^{4}\,{b_{1}}\,{a_{0}} \nonumber \\
\mathrm{T_{7}} = 2\,{b_{0}} + 12\,{b_{1}}\,{a_{0}}\,H^{2} \nonumber \\
\mathrm{T_{8}} = 2\,H\,{b_{0}} + 4\,H^{3}\,{b_{1}}\,{a_{0}}
\nonumber \\
\mathrm{T_{9}} = 48\,{b_{1}}\,{a_{0}}\,H\,\mathrm{cos}(H) -
6\,\mathrm{T_{7}} \,\mathrm{sin}(H) \label{nm3c}
\end{eqnarray}

We demand that the phase-lag and its first, second and third derivative
are equal to zero and we find out that:

\begin{eqnarray}
{a_{0}} = {\displaystyle \frac {1}{4}} \Bigl (
3\,\mathrm{cos}(H)^{2} + \mathrm{cos}(H)^{2}\,H^{2} + 2\,H^{2} - 3
\Bigr ) \nonumber \\ / \Bigl ( \Bigl ( 6\,\mathrm{
cos}(H)^{3}\,H + 6\,\mathrm{sin}(H)\,\mathrm{cos}(H)^{2} \nonumber \\
 - 2\,\mathrm{cos}(H)^{2}\,H^{2}\,\mathrm{sin}(H) +
\mathrm{cos}(H)^{2}\,H^{3} \nonumber \\
+ 3\,\mathrm{cos}(H)^{2}\,H - 6\,
\mathrm{cos}(H)\,\mathrm{sin}(H) \nonumber \\
 - 4\,\mathrm{cos}(H)\,H^{2}\,\mathrm{sin}(H) - 12\,
\mathrm{cos}(H)\,H + 2\,H^{3} \nonumber \\ + 3\,H +
12\,\mathrm{sin}(H)\,H^{2}
\Bigr ) H \Bigr ) \nonumber \\
{c_{1}} =  - 2 \, \Bigl ( - 12\,\mathrm{cos}(H)^{3}\,H +
\mathrm{cos}(H)^{ 2}\,H^{3} - 21\,\mathrm{cos}(H)^{2}\,H \nonumber
\\ - 12\,\mathrm{sin}(H)\, \mathrm{cos}(H)^{2} -
4\,\mathrm{cos}(H)^{2}\,H^{2}\,\mathrm{sin}(H) + 12\,
\mathrm{cos}(H)\,\mathrm{sin}(H) \nonumber \\
- 8\,\mathrm{cos}(H)\,H^{2}\,\mathrm{sin}(H) +
24\,\mathrm{cos}(H)\,H + 2\,H^{3} + 9\,H +
24\,\mathrm{sin}(H)\,H^{2} \Bigr ) / \nonumber \\
\Bigl ( \mathrm{cos}(H)^{2}\,H^{3} - 21\,\mathrm{cos}(H)^{2}\,H +
8\, \mathrm{cos}(H)\,H^{2}\,\mathrm{sin}(H) \nonumber \\
 - 12\,\mathrm{cos}(H)\,H - 12\,\mathrm{cos}(H)\,\mathrm{sin}(H) + 4\,\mathrm{sin}(H)\,H^{2} \nonumber \\
 + 33\,H + 12\,\mathrm{sin}(H) + 2\,H^{3} \Bigr ) \nonumber \\
{b_{0}} =  - 2 \, \Bigl ( 3\,\mathrm{cos}(H)^{2}\,H +
\mathrm{cos}(H)^{2}\, H^{3} + 6\,\mathrm{cos}(H)\,\mathrm{sin}(H)
\nonumber \\
+ 4\,\mathrm{cos}(H) \,H^{2}\,\mathrm{sin}(H) + 6\,\mathrm{cos}(H)\,H \nonumber \\
 + 2\,\mathrm{sin}(H)\,H^{2} - 9\,H - 6\,\mathrm{sin}(H)
 + 2\,H^{3} \Bigr ) / \nonumber \\
 \Bigl ( \Bigl ( \mathrm{cos}(H)^{2}\,H^{3} - 21\,\mathrm{cos}(H)^{2}\,H  + 8\,\mathrm{cos}(H)\,H^{2}\,\mathrm{sin}(H) \nonumber \\
- 12\, \mathrm{cos}(H)\,H - 12\,\mathrm{cos}(H)\,\mathrm{sin}(H)
\nonumber \\
+ 4\,\mathrm{sin}(H)\,H^{2} + 33\,H + 12\,\mathrm{sin}(H) + 2\,H^{3} \Bigr ) H^{2} \Bigr ) \nonumber \\
{b_{1}} = 4\, \Bigl ( 6\,\mathrm{cos}(H)^{3}\,H +
6\,\mathrm{sin}(H)\, \mathrm{cos}(H)^{2} -
2\,\mathrm{cos}(H)^{2}\,H^{2}\,\mathrm{sin}
(H) \nonumber \\
+ \mathrm{cos}(H)^{2}\,H^{3} + 3\,\mathrm{cos}(H)^{2}\,H -
6\,\mathrm{cos}(H)\,\mathrm{sin}(H) \nonumber \\
- 4\,\mathrm{cos}(H)\,H^{2}\,\mathrm{sin}(H) -
12\,\mathrm{cos}(H)\,H + 2\,H^{3} + 3
\,H + 12\,\mathrm{sin}(H)\,H^{2} \Bigr ) / \nonumber \\
 \Bigl ( \Bigl ( \mathrm{cos}(H)^{2}\,H^{3} - 21\,\mathrm{cos}(H)^{2}\,H + 8
\,\mathrm{cos}(H)\,H^{2}\,\mathrm{sin}(H) - 12\,\mathrm{cos}(H)\,
H \nonumber \\
 - 12\,\mathrm{cos}(H)\,\mathrm{sin}(H) + 4\,\mathrm{sin}(H)\,H^{2} + 33\,H + 12\,\mathrm{sin}(H) + 2\,H^{3} \Bigr ) H^{2}
\Bigr ) \label{nm4c}
\end{eqnarray}

For small values of $|H|$ the formulae given by (\ref{nm4c}) are
subject to heavy cancellations. In this case the following Taylor
series expansions should be used:

\begin{eqnarray}
a_{0} = {\displaystyle \frac {1}{200}}  + {\displaystyle \frac
{1}{1260}} \,H^{2} + {\displaystyle \frac {29}{504000}} \,H ^{4} +
{\displaystyle \frac {1433}{1164240000}} \,H^{6} \nonumber \\ -
{\displaystyle \frac {63101}{363242880000}} \,H^{8}  -
{\displaystyle \frac {2228861}{127135008000000}} \,H^{ 10} \nonumber
\\- {\displaystyle \frac {8804897}{77806624896000000}} \,H^{12}
 + {\displaystyle \frac {240953700049}{
2048959660011264000000}} \,H^{14} \nonumber \\ + {\displaystyle
\frac { 9699610781879}{819583864004505600000000}} \,H^{16} + \ldots \nonumber \\
c_{1} = - 2 + {\displaystyle \frac {1}{6048}} \,H^{8} +
{\displaystyle \frac {1}{43200}} \,H^{10} + {\displaystyle \frac
{1}{532224}} \,H^{12} \nonumber \\ + {\displaystyle \frac {41}{
5943974400}} \,H^{14} - {\displaystyle \frac {601}{24141680640}}
\,H^{16} + \ldots \nonumber \\
b_{0} = {\displaystyle \frac {1}{12}}  - {\displaystyle \frac
{1}{1008}} \,H^{4} - {\displaystyle \frac {31}{181440}} \,H ^{6}
\nonumber \\  - {\displaystyle \frac {221}{13685760}} \,H^{8} -
{\displaystyle \frac {619}{1345344000}} \,H^{10} \nonumber \\ +
{\displaystyle \frac {25031}{174356582400}} \,H^{12}
 + {\displaystyle \frac {84256583}{2667655710720000}} \,H^{14} \nonumber \\
 + {\displaystyle \frac {1030007057}{290289444157440000}} \,
H^{16} + \ldots \nonumber \\
b_{1} = {\displaystyle \frac {5}{6}}  + {\displaystyle \frac
{1}{504}} \,H^{4} - {\displaystyle \frac {29}{90720}} \,H^{6}
\nonumber \\
- {\displaystyle \frac {3271}{47900160}} \,H^{8} - {\displaystyle
\frac {35293}{4540536000}} \,H^{10} \nonumber \\
- {\displaystyle \frac {36019}{87178291200}} \,H^{12}
 + {\displaystyle \frac {47333617}{1333827855360000}} \,H^{14} \nonumber \\
 + {\displaystyle \frac {294008389}{24562952967168000}} \,H^{16}
 + \ldots \label{nm4tc}
\end{eqnarray}

The behavior of the coefficients is given in the following Figure 4.

\setcounter{figure}{3}
\begin{figure}[htb!]
\begin{center}
\begin{tabular}{ccc}
\includegraphics[width=0.4\textwidth]{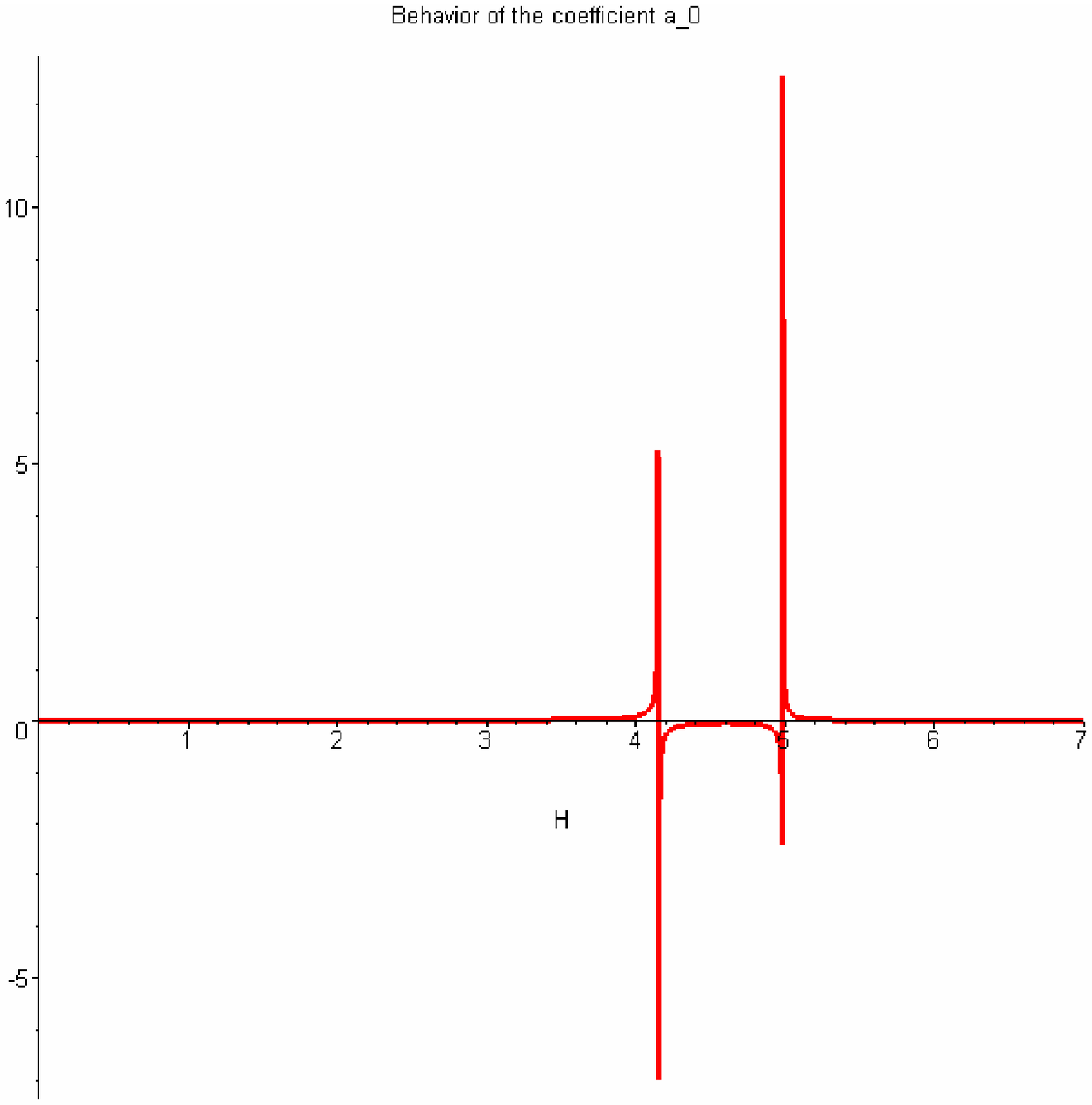}&\quad &
\includegraphics[width=0.4\textwidth]{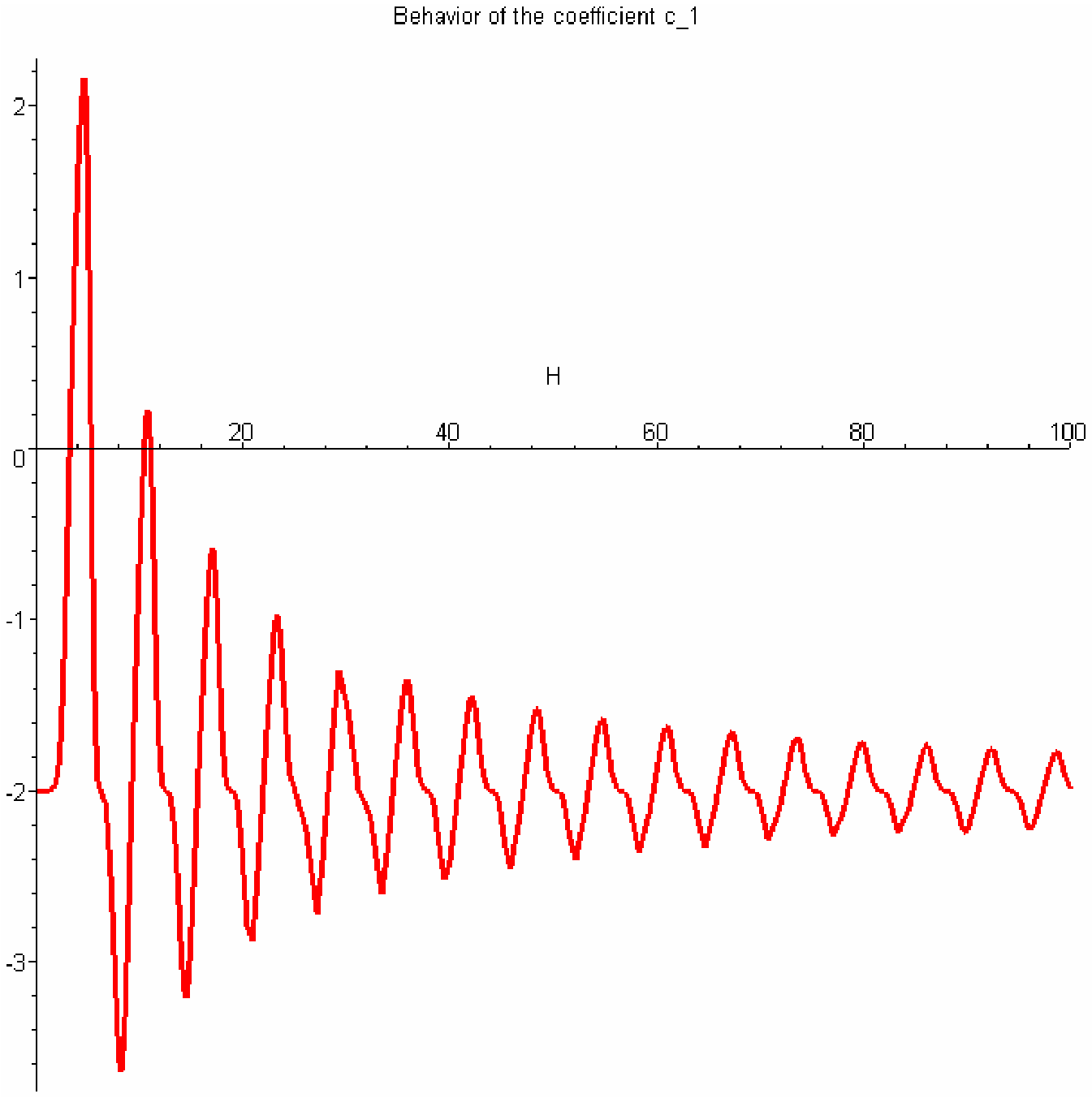}\\
\end{tabular}
\begin{tabular}{ccc}
\includegraphics[width=0.4\textwidth]{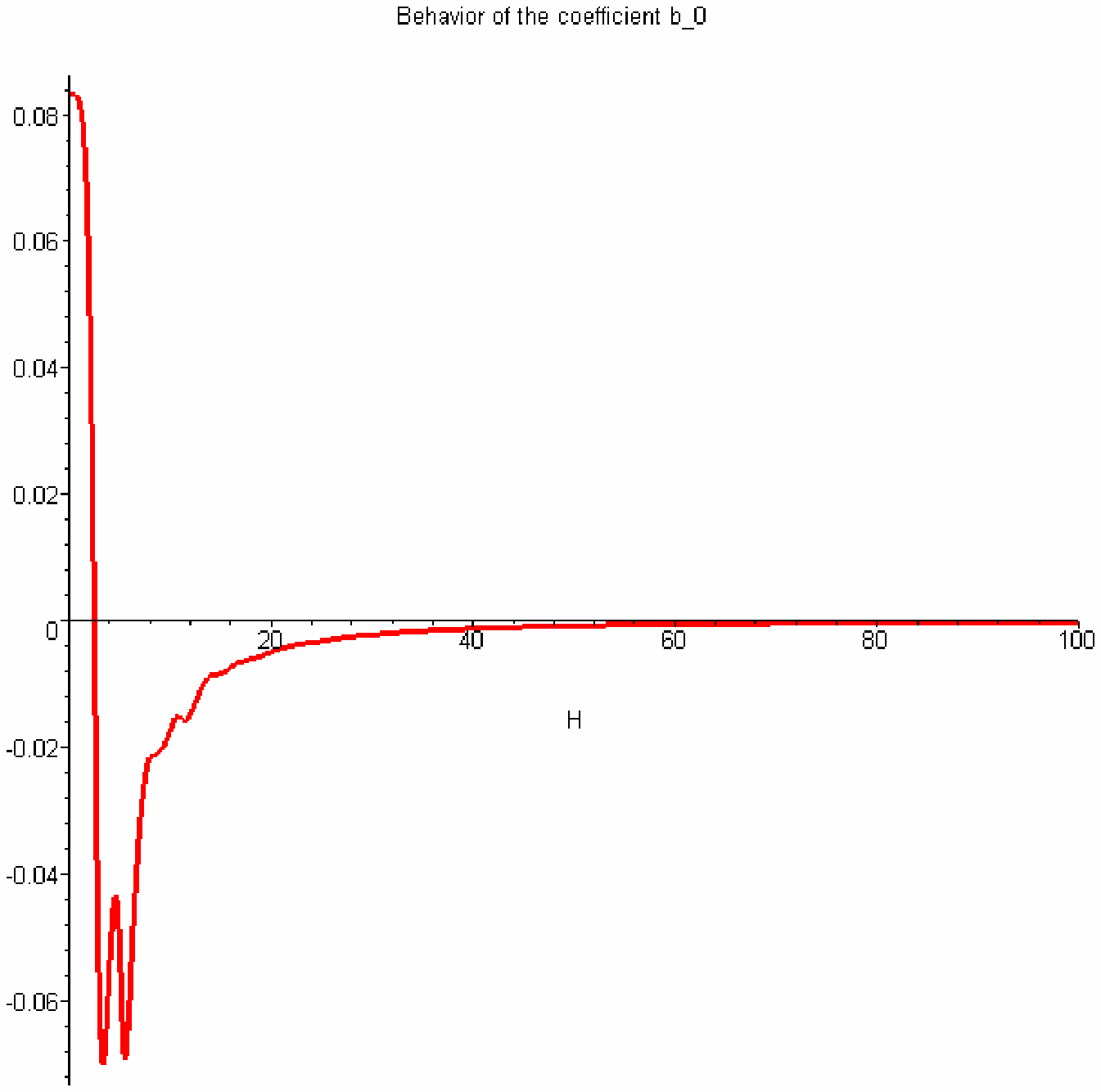}&\quad &
\includegraphics[width=0.4\textwidth]{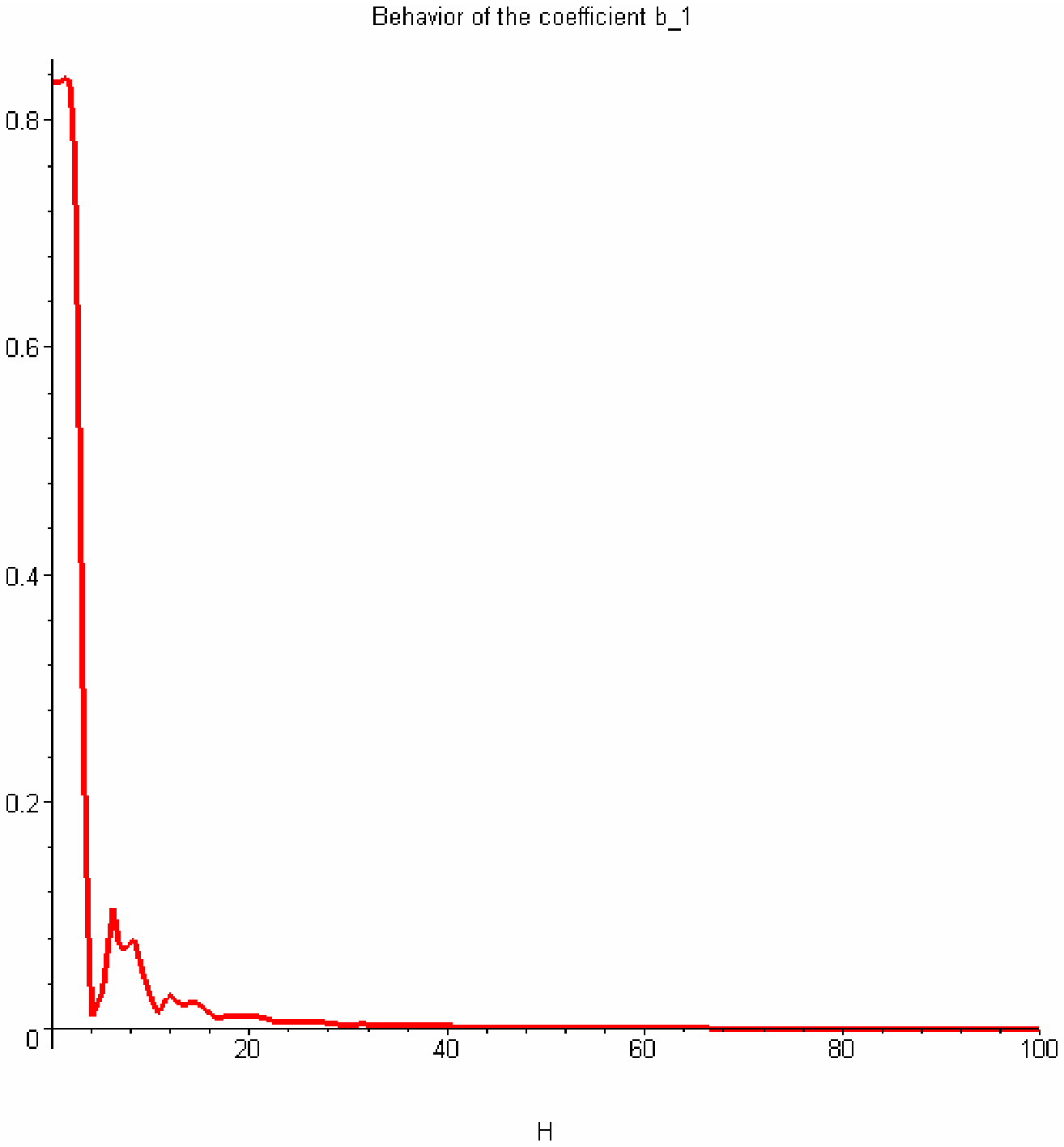}\\
\end{tabular}
\caption[]{\label{fig4}
Behavior of the coefficients of the new method given by (\ref{nm4c})
for several values of $H$.} \normalsize
\end{center}
\end{figure}

\vspace*{0.7cm}

The local truncation error of the new proposed method is given by:

\begin{equation}
\mathrm{LTE} = \frac{h^8}{6048} \left (y_{n}^{(8)} + 4\,
\omega^{2}\,\,y_{n}^{(6)} + 6\, \omega^{4}\,\,y_{n}^{(4)} + 4
\omega^6 \, y_{n}^{(2)} + \omega^8 \, y_{n} \right ) \label{ltec}
\end{equation}

\section{Error Analysis}

We will study the following methods:

\begin{itemize}
\item The First Method of the Family (mentioned as $PL1$)  \item The Second Method of the Family (mentioned as $PL2$) \item The Third Method of the Family (mentioned as $PL3$) \item The Fourth Method of the
Family (mentioned as $PL4$)
\end{itemize}
The error analysis is based on the following steps:

\begin{itemize}

\item The radial time independent Schr\"odinger equation is of the
form

\begin{equation}
y''(x) = f(x) \, y(x) \label{ierr1}
\end{equation}

\item Based on the paper of Ixaru and Rizea \cite{ixaru85}, the
function f(x) can be written in the form:

\begin{equation}
f(x) = g(x) + G \label{err2}
\end{equation}
\noindent where $g(x) = V(x) - V_{c} = g$, where $V_{c}$ is the
constant approximation of the potential and $G = v^2 = V_{c} - E$.

\item We express the derivatives $y_{n}^{(i)}, \, i=2,3,4, \ldots ,$
which are terms of the local truncation error formulae, in terms
of the equation (\ref{err1}). The expressions are presented as
polynomials of $G$.

\item Finally, we substitute the expressions of the derivatives,
produced in the previous step, into the local truncation error
formulae.

\end{itemize}

Based on the procedure mentioned above and on the formulae:

\begin{eqnarray*}
y_{n}^{(2)} = (\mathrm{V}(x) - \mathit{V_{c}} + G)\,\mathrm{y}(x) \nonumber \\
y_{n}^{(4)} = ({\frac {d^{2}}{dx^{2}}}\,\mathrm{V}(x))\,
\mathrm{y}(x) + 2\,({\frac {d}{dx}}\,\mathrm{V}(x))\,({\frac {d}{
dx}}\,\mathrm{y}(x)) \nonumber \\ + (\mathrm{V}(x) -
\mathit{V_{c}} +
G)\,( {\frac {d^{2}}{dx^{2}}}\,\mathrm{y}(x)) \nonumber \\
y_{n}^{(6)} = ({\frac {d^{4}}{dx^{4}}}\,\mathrm{V}(x))\,
\mathrm{y}(x) + 4\,({\frac {d^{3}}{dx^{3}}}\,\mathrm{V}(x))\,(
{\frac {d}{dx}}\,\mathrm{y}(x)) \nonumber \\ + 3\,({\frac
{d^{2}}{dx^{2}}}\, \mathrm{V}(x))\,({\frac
{d^{2}}{dx^{2}}}\,\mathrm{y}(x)) \nonumber
\\ + 4\,( {\frac {d}{dx}}\,\mathrm{V}(x))^{2}\,\mathrm{y}(x)
\nonumber \\ + 6\,(\mathrm{V}(x) - \mathit{V_{c}} + G)\,({\frac
{d}{dx}} \,\mathrm{y}(x))\,({\frac {d}{dx}}\,\mathrm{V}(x))
\nonumber \\ + 4\,(\mathrm{ U}(x) - \mathit{V_{c}} +
G)\,\mathrm{y}(x)\,({\frac {d^{2}}{dx^{2}}} \,\mathrm{V}(x))
\nonumber \\ + (\mathrm{V}(x) - \mathit{V_{c}} + G)^{2}\,({\frac
{d^{2}}{ dx^{2}}}\,\mathrm{y}(x)) \ldots
\end{eqnarray*}
\noindent we obtain the following expressions:

\subsection*{The First Method of the Family}

\begin{eqnarray}
\mathrm{LTE_{PL1}} =h^8 \, \Bigl [ - {\displaystyle \frac {1}{6048}}
\,\mathrm{g}(x)\,\mathrm{y}(x)\,G^{3}  + \Bigl ( - {\displaystyle
\frac {5}{2016}} \,({\frac {d^{2}}{d
x^{2}}}\,\mathrm{g}(x))\,\mathrm{y}(x) \nonumber \\ - {\displaystyle
\frac {1 }{1008}} \,({\frac {d}{dx}}\,\mathrm{g}(x))\,({\frac
{d}{dx}}\, \mathrm{y}(x)) - {\displaystyle \frac
{1}{2016}} \,\mathrm{g}(x)^{2}\,\mathrm{y}(x) \Bigr ) \,G^{2} + \Bigl ( \nonumber \\
 - {\displaystyle \frac {5}{2016}} \,({\frac {d^{4}}{dx^{4}}}\,
\mathrm{g}(x))\,\mathrm{y}(x) - {\displaystyle \frac {5}{1512}}
\,({\frac {d^{3}}{dx^{3}}}\,\mathrm{g}(x))\,({\frac {d}{dx}}\,
\mathrm{y}(x)) \nonumber \\ - {\displaystyle \frac {1}{336}}
\,\mathrm{g}(x)\, ({\frac {d}{dx}}\,\mathrm{y}(x))\,({\frac
{d}{dx}}\,\mathrm{g}(x)) - {\displaystyle \frac {37}{6048}}
\,\mathrm{g}(x)\, \mathrm{y}(x)\,({\frac
{d^{2}}{dx^{2}}}\,\mathrm{g}(x)) \nonumber
\\ - {\displaystyle \frac {1}{252}} \,({\frac
{d}{dx}}\,\mathrm{g}(x)) ^{2}\,\mathrm{y}(x) - {\displaystyle \frac
{1}{2016}} \,\mathrm{g}(x)^{3}\,\mathrm{y}(x) \Bigr ) \, G \nonumber \\
 - {\displaystyle \frac {1}{6048}} \,({\frac {d^{6}}{dx^{6
}}}\,\mathrm{g}(x))\,\mathrm{y}(x) - {\displaystyle \frac {1}{
1008}} \,({\frac{d^{5}}{dx^{5}}}\,\mathrm{g}(x))\,({\frac {d}{dx
}}\,\mathrm{y}(x))\nonumber \\ - {\displaystyle \frac
{1}{378}}\,\mathrm{g}(x)\,\mathrm{y}(x)\,({\frac
{d^{4}}{dx^{4}}}\,\mathrm{g}(x)) - {\displaystyle \frac {5}{2016}}
\,({\frac
{d^{2}}{dx^{2}}}\,\mathrm{g}(x))^{2}\,\mathrm{y}(x)\nonumber \\ -
{\displaystyle \frac {13 }{3024}} \,({\frac
{d}{dx}}\,\mathrm{g}(x))\,\mathrm{y}(x)\,( {\frac
{d^{3}}{dx^{3}}}\,\mathrm{g}(x))  - {\displaystyle \frac {1}{252}}
\,\mathrm{g}(x)\,( {\frac{d}{dx}}\,\mathrm{y}(x))\,({\frac {d^{3}}{dx^{3}}}\, \mathrm{g}(x)) \nonumber \\
- {\displaystyle \frac {1}{504}} \,\mathrm{g}(x)^{ 2}\,({\frac
{d}{dx}}\,\mathrm{y}(x))\,({\frac {d}{dx}}\,\mathrm{g }(x))
 \nonumber \\ - {\displaystyle \frac {1}{126}} \,({\frac {d}{dx}}\,
\mathrm{g}(x))\,({\frac {d}{dx}}\,\mathrm{y}(x))\,({\frac {d^{2}
}{dx^{2}}}\,\mathrm{g}(x))\nonumber \\ - {\displaystyle \frac
{11}{3024}} \, \mathrm{g}(x)^{2}\,\mathrm{y}(x)\,({\frac
{d^{2}}{dx^{2}}}\,\mathrm{g}(x))  - {\displaystyle \frac {1}{216}}
\,\mathrm{g}(x)\, \mathrm{y}(x)\,({\frac
{d}{dx}}\,\mathrm{g}(x))^{2} \nonumber \\- {\displaystyle \frac
{1}{6048}} \,\mathrm{g}(x)^{4}\,\mathrm{y}(x) \Bigr ]  \label{err3a}
\end{eqnarray}

\subsection*{The Second Method of the Family}

\begin{eqnarray}
\mathrm{LTE_{PL2}} = h^8 \, \Bigl [ ({\displaystyle \frac {19}{9072
}} \,({\frac {d^{2}}{dx^{2}}}\,\mathrm{g}(x))\,\mathrm{y}(x) +
{\displaystyle \frac {1}{1512}} \,({\frac {d}{dx}}\,\mathrm{g}(x)
)\,({\frac {d}{dx}}\,\mathrm{y}(x)) \nonumber \\ + {\displaystyle
\frac {1}{ 3024}} \,\mathrm{g}(x)^{2}\,\mathrm{y}(x))\,G^{2} + (
{\displaystyle \frac {11}{4536}} \,({\frac {d^{4}}{dx^{4}}}\,
\mathrm{g}(x))\,\mathrm{y}(x) \nonumber \\
+ {\displaystyle \frac{1}{324}} \, ({\frac
{d^{3}}{dx^{3}}}\,\mathrm{g}(x))\,({\frac {d}{dx}}\, \mathrm{y}(x))
+ {\displaystyle \frac {1}{378}} \,\mathrm{g}(x)\, ({\frac
{d}{dx}}\,\mathrm{y}(x))\,({\frac {d}{dx}}\,\mathrm{g}(x) )
\nonumber \\ + {\displaystyle \frac {13}{2268}} \,\mathrm{g}(x)\,
\mathrm{y}(x)\,({\frac {d^{2}}{dx^{2}}}\,\mathrm{g}(x)) +
{\displaystyle \frac {17}{4536}} \,({\frac {d}{dx}}\,\mathrm{g}(x
))^{2}\,\mathrm{y}(x) \nonumber \\ + {\displaystyle \frac {1}{2268}}
\,\mathrm{g}(x)^{3}\,\mathrm{y}(x))G  + {\displaystyle \frac
{1}{6048}} \,({\frac {d^{6}}{dx^{6}}}\,\mathrm{g}(x))\,\mathrm{y}(x) \nonumber \\
+ {\displaystyle \frac {1}{ 1008}} \,({\frac
{d^{5}}{dx^{5}}}\,\mathrm{g}(x))\,({\frac {d}{dx }}\,\mathrm{y}(x))
+ {\displaystyle \frac {1}{378}} \,\mathrm{g}(
x)\,\mathrm{y}(x)\,({\frac {d^{4}}{dx^{4}}}\,\mathrm{g}(x)) \nonumber \\
 + {\displaystyle \frac {5}{2016}} \,({\frac {d^{2}}{dx^{2
}}}\,\mathrm{g}(x))^{2}\,\mathrm{y}(x) + {\displaystyle \frac
{13}{3024}} \,({\frac {d}{dx}}\,\mathrm{g}(x))\,\mathrm{y}(x)\,(
{\frac {d^{3}}{dx^{3}}}\,\mathrm{g}(x)) \nonumber \\
 + {\displaystyle \frac {1}{252}} \,\mathrm{g}(x)\,( {\frac
{d}{dx}}\,\mathrm{y}(x))\,({\frac {d^{3}}{dx^{3}}}\, \mathrm{g}(x))
+ {\displaystyle \frac {1}{504}} \,\mathrm{g}(x)^{ 2}\,({\frac
{d}{dx}}\,\mathrm{y}(x))\,({\frac {d}{dx}}\,\mathrm{g
}(x)) \nonumber \\
 + {\displaystyle \frac {1}{126}} \,({\frac {d}{dx}}\,
\mathrm{g}(x))\,({\frac {d}{dx}}\,\mathrm{y}(x))\,({\frac {d^{2}
}{dx^{2}}}\,\mathrm{g}(x)) + {\displaystyle \frac {11}{3024}} \,
\mathrm{g}(x)^{2}\,\mathrm{y}(x)\,({\frac {d^{2}}{dx^{2}}}\,
\mathrm{g}(x)) \nonumber \\
 + {\displaystyle \frac {1}{216}} \,\mathrm{g}(x)\,
\mathrm{y}(x)\,({\frac {d}{dx}}\,\mathrm{g}(x))^{2} + {\displaystyle
\frac {1}{6048}} \,\mathrm{g}(x)^{4}\,\mathrm{y}(x) \Bigr ]
\label{err4a}
\end{eqnarray}

\subsection*{The Third Method of the Family}

\begin{eqnarray}
\mathrm{LTE_{PL3}} = h^8 \, \Bigl [ {\displaystyle \frac {1}{756}}
\,({\frac
{d^{2}}{dx^{2}}}\,\mathrm{g}(x))\,\mathrm{y}(x)\,G^{2} \nonumber \\
+ \Bigl ( {\displaystyle \frac {1}{432}} \,({\frac
{d^{4}}{dx^{4}}}\, \mathrm{g}(x))\,\mathrm{y}(x) + {\displaystyle
\frac {1}{378}} \, ({\frac {d^{3}}{dx^{3}}}\,\mathrm{g}(x))\,({\frac
{d}{dx}}\, \mathrm{y}(x)) \nonumber \\
 + {\displaystyle \frac {1}{504}} \,\mathrm{g}(x)\,( {\frac
{d}{dx}}\,\mathrm{y}(x))\,({\frac {d}{dx}}\,\mathrm{g}(x)) \nonumber
\\ + {\displaystyle \frac {5}{1008}} \,\mathrm{g}(x)\,\mathrm{y}(x)
\,({\frac {d^{2}}{dx^{2}}}\,\mathrm{g}(x)) + {\displaystyle \frac
{5}{1512}} \,({\frac {d}{dx}}\,\mathrm{g}(x))^{2}\,\mathrm{ y}(x)
\nonumber \\
 + {\displaystyle \frac {1}{3024}} \,\mathrm{g}(x)^{3}\,
\mathrm{y}(x) \Bigr ) \, G + {\displaystyle \frac {1}{6048}} \,(
{\frac {d^{6}}{dx^{6}}}\,\mathrm{g}(x))\,\mathrm{y}(x) \nonumber \\
+ {\displaystyle \frac {1}{1008}} \,({\frac {d^{5}}{dx^{5}}}\,
\mathrm{g}(x))\,({\frac {d}{dx}}\,\mathrm{y}(x)) \nonumber \\
 + {\displaystyle \frac {1}{378}} \,\mathrm{g}(x)\,
\mathrm{y}(x)\,({\frac {d^{4}}{dx^{4}}}\,\mathrm{g}(x)) \nonumber
\\+ {\displaystyle \frac {5}{2016}} \,({\frac {d^{2}}{dx^{2}}}\,
\mathrm{g}(x))^{2}\,\mathrm{y}(x) + {\displaystyle \frac {13}{
3024}} \,({\frac {d}{dx}}\,\mathrm{g}(x))\,\mathrm{y}(x)\,( {\frac
{d^{3}}{dx^{3}}}\,\mathrm{g}(x)) \nonumber \\
 + {\displaystyle \frac {1}{252}} \,\mathrm{g}(x)\,( {\frac
{d}{dx}}\,\mathrm{y}(x))\,({\frac {d^{3}}{dx^{3}}}\, \mathrm{g}(x)) \nonumber \\
+ {\displaystyle \frac {1}{504}} \,\mathrm{g}(x)^{ 2}\,({\frac
{d}{dx}}\,\mathrm{y}(x))\,({\frac {d}{dx}}\,\mathrm{g }(x))
\nonumber \\
 + {\displaystyle \frac {1}{126}} \,({\frac {d}{dx}}\,
\mathrm{g}(x))\,({\frac {d}{dx}}\,\mathrm{y}(x))\,({\frac {d^{2}
}{dx^{2}}}\,\mathrm{g}(x)) + {\displaystyle \frac {11}{3024}} \,
\mathrm{g}(x)^{2}\,\mathrm{y}(x)\,({\frac {d^{2}}{dx^{2}}}\,
\mathrm{g}(x)) \nonumber \\
 + {\displaystyle \frac {1}{216}} \,\mathrm{g}(x)\,
\mathrm{y}(x)\,({\frac {d}{dx}}\,\mathrm{g}(x))^{2} + {\displaystyle
\frac {1}{6048}} \,\mathrm{g}(x)^{4}\,\mathrm{y}(x) \Bigr ]
\label{err5b}
\end{eqnarray}

\subsection*{The Fourth Method of the Family}

\begin{eqnarray}
\mathrm{LTE_{PL4}} =h^8 \, \Bigl [ \Bigl ( {\displaystyle \frac
{1}{504} } \,({\frac {d^{4}}{dx^{4}}}\,\mathrm{g}(x))\,\mathrm{y}(x)
+ {\displaystyle \frac {1}{756}} \,({\frac {d^{3}}{dx^{3}}}\,
\mathrm{g}(x))\,({\frac {d}{dx}}\,\mathrm{y}(x)) \nonumber \\  +
{\displaystyle \frac {1}{378}} \,\mathrm{g}(x)\,\mathrm{y}(x)\,(
{\frac {d^{2}}{dx^{2}}}\,\mathrm{g}(x)) + {\displaystyle \frac
{1}{504}} \,({\frac {d}{dx}}\, \mathrm{g}(x))^{2}\,\mathrm{y}(x)
\Bigr ) \, G \nonumber \\  + {\displaystyle
\frac{1}{6048}}\,({\frac{d^{6}}{dx^{6}}}\,\mathrm{g}(x))\,
\mathrm{y}(x) + {\displaystyle \frac {1}{1008}} \,({\frac {d^{5}
}{dx^{5}}}\,\mathrm{g}(x))\,({\frac {d}{dx}}\,\mathrm{y}(x))
\nonumber \\  + {\displaystyle \frac {1}{378}} \,\mathrm{g}(x)\,
\mathrm{y}(x)\,({\frac {d^{4}}{dx^{4}}}\,\mathrm{g}(x)) +
{\displaystyle \frac {5}{2016}} \,({\frac {d^{2}}{dx^{2}}}\,
\mathrm{g}(x))^{2}\,\mathrm{y}(x) \nonumber \\  + {\displaystyle
\frac {13}{ 3024}} \,({\frac
{d}{dx}}\,\mathrm{g}(x))\,\mathrm{y}(x)\,( {\frac
{d^{3}}{dx^{3}}}\,\mathrm{g}(x)) \nonumber \\  + {\displaystyle
\frac {1}{252}} \,\mathrm{g}(x)\,( {\frac
{d}{dx}}\,\mathrm{y}(x))\,({\frac {d^{3}}{dx^{3}}}\, \mathrm{g}(x))
+ {\displaystyle \frac {1}{504}} \,\mathrm{g}(x)^{ 2}\,({\frac
{d}{dx}}\,\mathrm{y}(x))\,({\frac {d}{dx}}\,\mathrm{g }(x))
\nonumber \\  + {\displaystyle \frac {1}{126}} \,({\frac {d}{dx}}\,
\mathrm{g}(x))\,({\frac {d}{dx}}\,\mathrm{y}(x))\,({\frac {d^{2}
}{dx^{2}}}\,\mathrm{g}(x)) \nonumber \\  + {\displaystyle \frac
{11}{3024}} \, \mathrm{g}(x)^{2}\,\mathrm{y}(x)\,({\frac
{d^{2}}{dx^{2}}}\, \mathrm{g}(x)) \nonumber \\  + {\displaystyle
\frac {1}{216}} \,\mathrm{g}(x)\, \mathrm{y}(x)\,({\frac
{d}{dx}}\,\mathrm{g}(x))^{2} + {\displaystyle \frac {1}{6048}}
\,\mathrm{g}(x)^{4}\,\mathrm{y}(x) \Bigr ]  \label{err5c}
\end{eqnarray}

We consider two cases in terms of the value of $E$:

\begin{itemize}

\item The Energy is close to the potential, i.e. $G = V_{c} - E
\approx 0$. So only the free terms of the polynomials in $G$ are
considered. Thus for these values of $G$, the methods are of
comparable accuracy. This is because the free terms of the
polynomials in $G$, are the same for the cases of the classical
method and of the new developed methods.

\item $G \gg 0$ or $G \ll 0$. Then $\mid G \mid$ is a large number. So,
we have the following asymptotic expansions of the equations
(\ref{err3a}), (\ref{err4a}), (\ref{err5b}) and (\ref{err5c}).

\subsection*{The First Method of the Family}

\begin{equation}
\mathrm{LTE_{PL1}} = h^8 \, \Bigl ( - {\displaystyle \frac
{1}{6048}} \,\mathrm{g}(x)\,\mathrm{y}(x)\,G^{3}  + \ldots \Bigr )
\label{err3b}
\end{equation}

\subsection*{The Second Method of the Family}

\begin{equation}
\begin{array}{l}
\mathrm{LTE_{PL2}} =h^8 \, \Bigl ({\displaystyle \frac {19}{9072 }}
\,({\frac {d^{2}}{dx^{2}}}\,\mathrm{g}(x))\,\mathrm{y}(x) +
{\displaystyle \frac {1}{1512}} \,({\frac {d}{dx}}\,\mathrm{g}(x)
)\,({\frac {d}{dx}}\,\mathrm{y}(x)) \nonumber \\
+ {\displaystyle
\frac {1}{ 3024}} \,\mathrm{g}(x)^{2}\,\mathrm{y}(x))\,G^{2} +
\ldots \Bigr ) \label{err4b}
\end{array}
\end{equation}

\subsection*{The Third Method of the Family}

\begin{eqnarray}
\mathrm{LTE_{PL3}} = h^8 \, \Bigl ( {\displaystyle \frac {1}{756}}
\,({\frac {d^{2}}{dx^{2}}}\,\mathrm{g}(x))\,\mathrm{y}(x)\,G^{2} +
\ldots \Bigr ) \label{err5d}
\end{eqnarray}

\subsection*{The Fourth Method of the Family}

\begin{eqnarray}
\mathrm{LTE_{PL4}} = h^8 \, \Bigl ( \Bigl ( {\displaystyle \frac
{1}{504} } \,({\frac {d^{4}}{dx^{4}}}\,\mathrm{g}(x))\,\mathrm{y}(x)
+ {\displaystyle \frac {1}{756}} \,({\frac {d^{3}}{dx^{3}}}\,
\mathrm{g}(x))\,({\frac {d}{dx}}\,\mathrm{y}(x)) \nonumber \\  +
{\displaystyle \frac {1}{378}} \,\mathrm{g}(x)\,\mathrm{y}(x)\,(
{\frac {d^{2}}{dx^{2}}}\,\mathrm{g}(x)) + {\displaystyle \frac
{1}{504}} \,({\frac {d}{dx}}\, \mathrm{g}(x))^{2}\,\mathrm{y}(x)
\Bigr ) \, G + \ldots \Bigr ) \label{err5e}
\end{eqnarray}

\end{itemize}

From the above equations we have the following theorem:

\begin{theorem}
For the First Method of the New Family of Methods the error
increases as the third power of $G$. For the Second and Third
Methods of the New Family of Methods the error increases as the
second power of $G$. For the Fourth Method of the New Family of
Methods the error increases as the first power of $G$. It is easy
one to see that the coefficient of the second power of $G$ in the
case of the second method of the New Family of Methods is
$1.583333333$ times larger than the coefficient of the second power
of $G$ in the case of the third method of the New Family of Methods.
So, for the numerical solution of the time independent radial
Schr\"odinger equation the new obtained Fourth Method of the New
Family of Methods is the most accurate one, especially for large
values of $\mid G \mid = \mid V_{c} - E \mid$.
\end{theorem}

\section{Stability Analysis}

We apply the new family of methods to the scalar test equation:

\begin{equation}
y'' = - t^2 y, \,
\end{equation}
\noindent where $t \neq \omega$. We obtain the following difference
equation:

\begin{equation}
\label{inm1s} \nonumber A_{1}(H,s) \, y_{n + 1} + A_{0}(H,s)\, y_{n}
+ A_{1}(H,s)\, y_{n - 1} = 0
\end{equation}

\noindent where $s = t\, h$, $h$ is the step length and $A_{0}(H,s)$
and $A_{1}(H,s)$ are polynomials of $s$.

The characteristic equation associated with (\ref{nm1s}) is given
by:

\begin{eqnarray}
\label{nm2s} A_{1}(H,s)\, s + A_{0}(H,s) + A_{1} (H,s)\,s^{ - 1} = 0
\end{eqnarray}

\noindent where

\begin{eqnarray}
{A_{1}(H,s)} = 1 + s^{2}\,{b_{0}} + s^{4}\,{b_{1}}\,{a_{0}} \nonumber \\
{A_{0}(H,s)} = {c_{1}} + s^{2}\,{b_{1}} - 2\,s^{4}\,{b_{1}}\,{a_{0}}
\label{nm6b}
\end{eqnarray}

\begin{defn} (see \cite{lambert})
A symmetric four-step method with the characteristic equation given
by (\ref{nm2s}) is said to have an {\it interval of periodicity}
$\left ( 0,w_{0}^{2} \right )$ if, for all $w \in \left (
0,w_{0}^{2} \right )$, the roots $z_{i}, \, i=1,2$ satisfy

\begin{equation}
z_{1,2} = e^{\pm i\,\theta(t\,h)}, \, | z_{i} | \le 1, \, i=3,4
\label{nm7}
\end{equation}

\noindent where $\theta(t\, h)$ is a real function of $t\,h$ and
$s=t\,h$ .
\end{defn}

\begin{defn} (see \cite{lambert})
A method is called P-stable if its interval of periodicity is
equal to $\left ( 0, \infty \right )$.
\end{defn}

\begin{theorem} (see \cite{rapsim91})
A symmetric two-step method with the characteristic equation given
by (\ref{nm2s}) is said to have a nonzero {\it interval of
periodicity} $\left ( 0,s_{0}^{2} \right )$ if, for all $s \in \left
( 0,s_{0}^{2} \right )$ the following relations are hold

\begin{equation}
P_{1}(H,s) > 0,\,  P_{2}(H,s) > 0, \label{nm7abc}
\end{equation}
\noindent where $H=\omega\,h$, $s=t\,h$ and:

\begin{eqnarray}
P_{1}(H,s) = A_{0}(H,s) + 2\, A_{1}(H,s) > 0, \nonumber \\
P_{2}(H,s) = A_{0}(H,s) - 2\, A_{1}(H,s) > 0, \label{nm7abc1}
\end{eqnarray}
\end{theorem}

\begin{defn}
A method is called singularly almost P-stable if its interval of
periodicity is equal to $(0,\infty) - S$\footnote{where $S$ is a set
of distinct points} only when the frequency of the phase fitting is
the same as the frequency of the scalar test equation, i.e. $H=s$.
\end{defn}

Based on (\ref{nm6b}) the stability polynomials (\ref{nm7abc1})
for the new developed methods take the form:

\begin{eqnarray}
P_{1}(H,s) =  {c_{1}} + v^{2}\,{b_{1}} + 2 + 2\,v^{2}\,{b_{0}}, \nonumber \\
P_{2}(H,s) = {c_{1}} + v^{2}\,{b_{1}} - 4\,v^{4}\,{b_{1}}\,{a_{0}} -
2 - 2\,v^{2}\,{b_{0}} \label{nm7abc2}
\end{eqnarray}

In Figures 5, 6, 7 and 8 we present the $s-H$ planes for the methods developed in this paper. A shadowed area denotes the $s-H$ region where the method is unstable, while a white area denotes the region where the method is stable. In Figure 5 we present the $s-H$ plane for the first method of the new family of method developed in this paper (paragraph 3.1). In Figure 6 we present the $s-H$ plane for the second method of the new family of method developed in this paper (paragraph 3.2). In Figure 7 we present the $s-H$ plane for the third method of the new family of method developed in this paper (paragraph 3.3). Finally, in Figure 8 we present the $s-H$ plane for the fourth method of the new family of method developed in this paper (paragraph 3.4).

\setcounter{figure}{4}
\begin{figure}[htb!]
\begin{center}
\includegraphics[width=\textwidth]{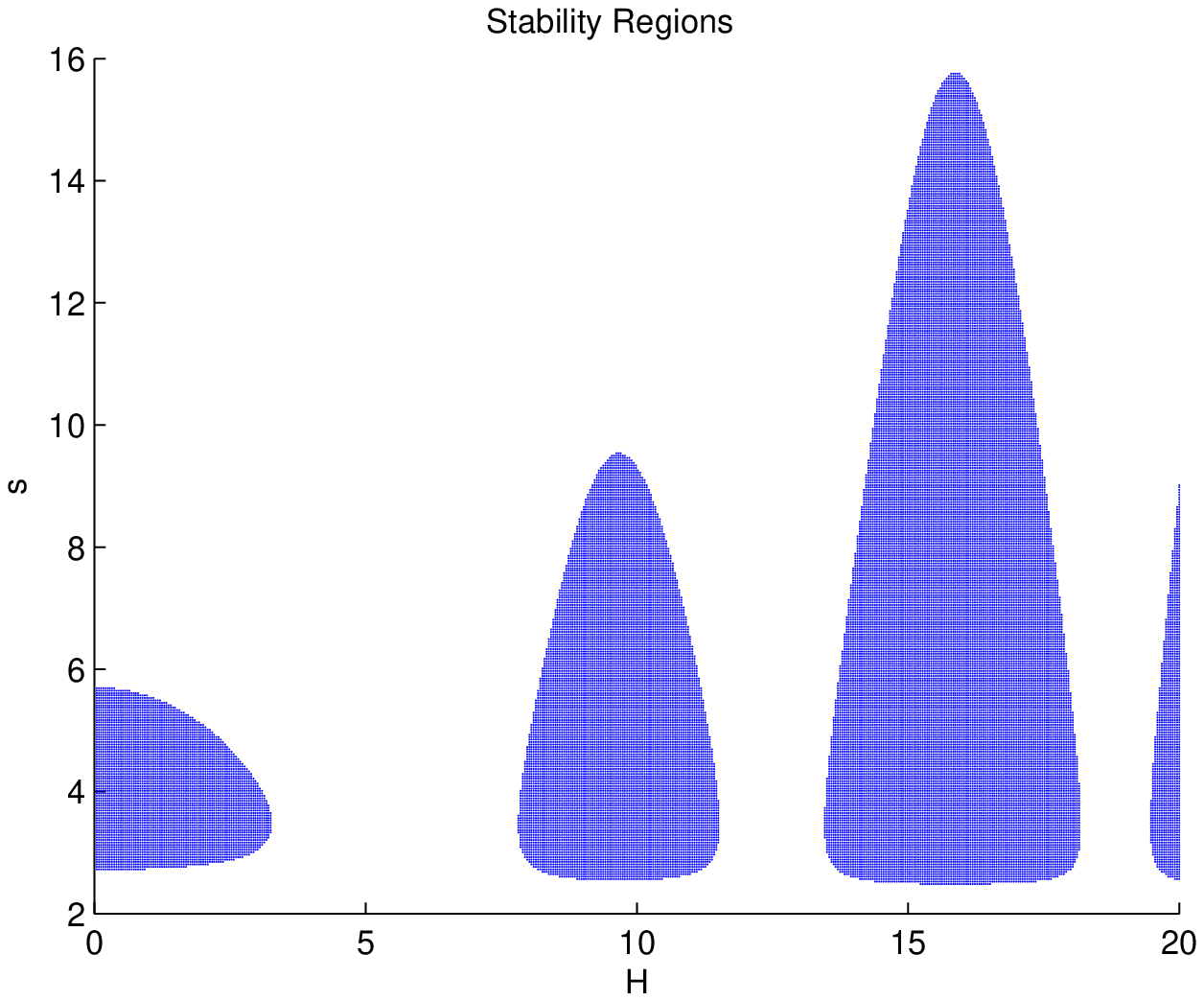}
\caption[]{\label{fig1s1}
$s-H$ plane of the first method of the new family of method
developed in this paper (paragraph 3.1)} \normalsize
\end{center}
\end{figure}

\vspace*{0.7cm}

\setcounter{figure}{5}
\begin{figure}[htb!]
\begin{center}
\includegraphics[width=\textwidth]{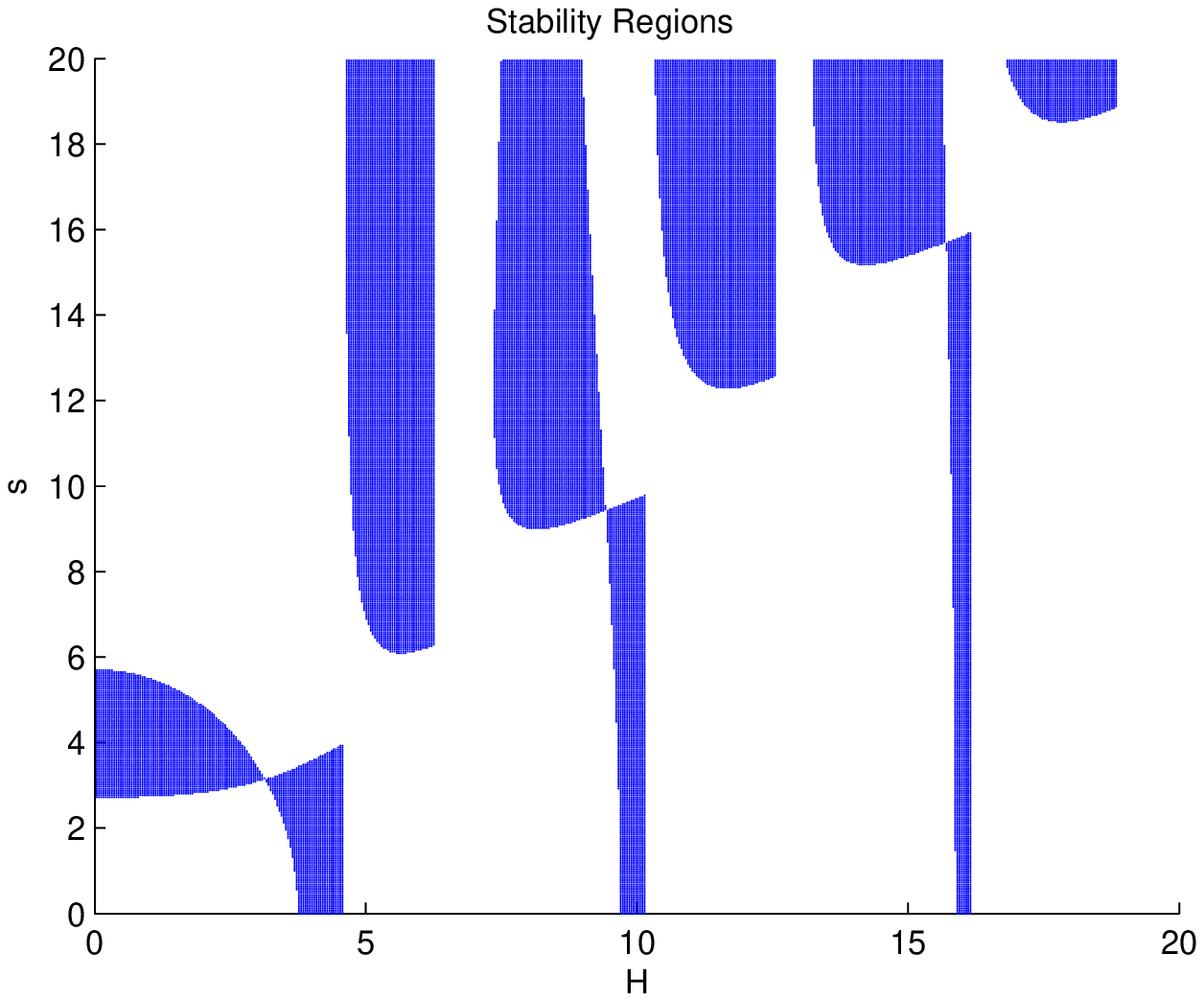}
\caption[]{\label{fig1s2}
$s-H$ plane of the second method of the new family of method
developed in this paper (paragraph 3.2)} \normalsize
\end{center}
\end{figure}

\vspace*{0.7cm}

\setcounter{figure}{6}
\begin{figure}[htb!]
\begin{center}
\includegraphics[width=\textwidth]{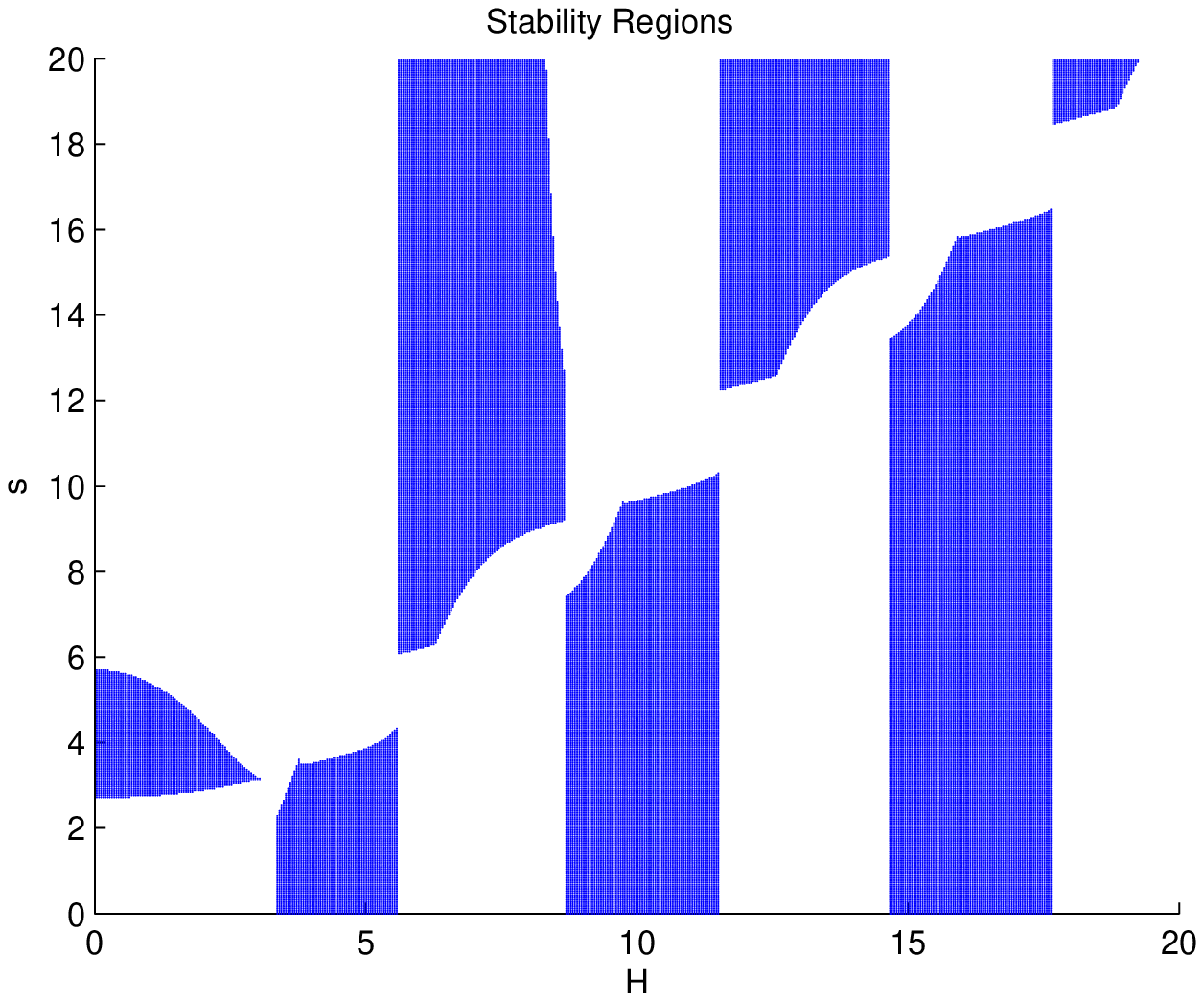}
\caption[]{\label{fig1s3}
$s-H$ plane of the third method of the new family of method
developed in this paper (paragraph 3.3)} \normalsize
\end{center}
\end{figure}

\vspace*{0.7cm}

\setcounter{figure}{7}
\begin{figure}[htb!]
\begin{center}
\includegraphics[width=0.9\textwidth]{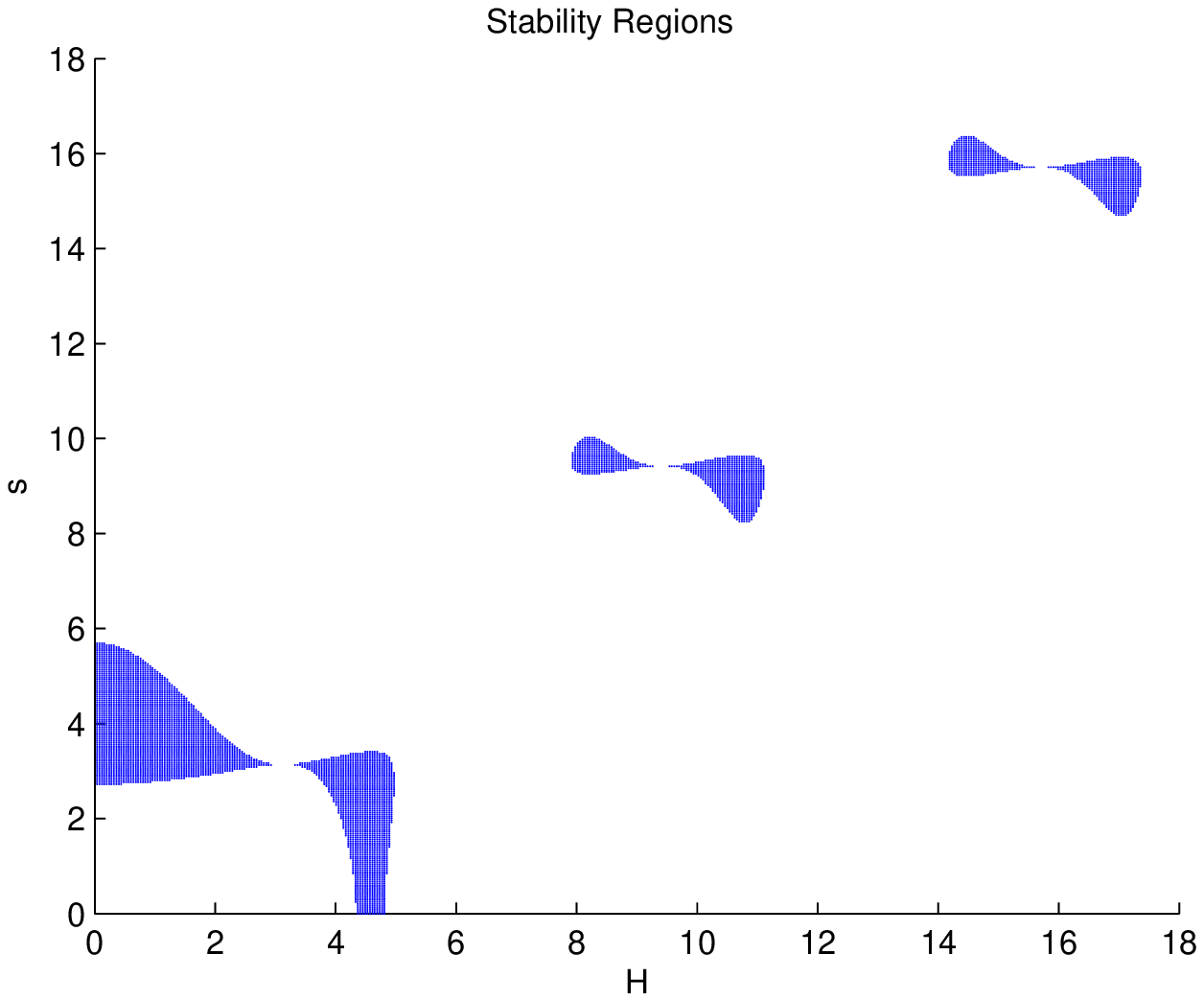}
\caption[]{\label{fig1s4}
$s-H$ plane of the fourth method of the new family of method
developed in this paper (paragraph 3.4)} \normalsize
\end{center}
\end{figure}

\vspace*{0.7cm}

\setcounter{figure}{8}
\begin{figure}[htb!]
\begin{center}
\begin{tabular}{ccc}
\includegraphics[width=0.4\textwidth]{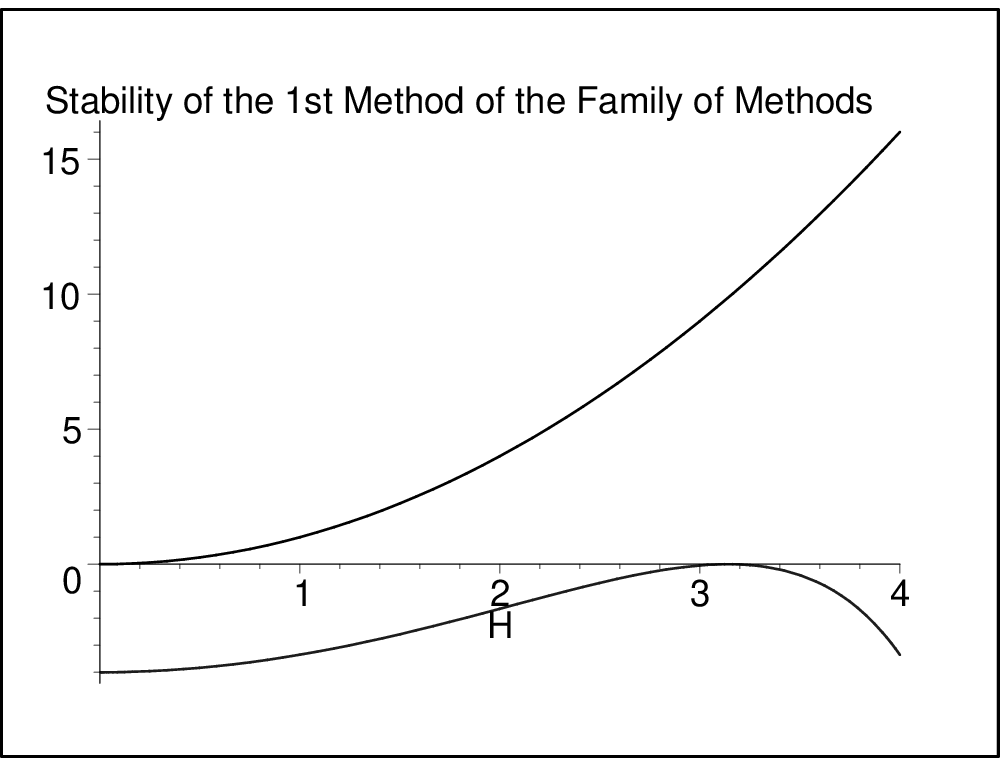}&\quad &
\includegraphics[width=0.4\textwidth]{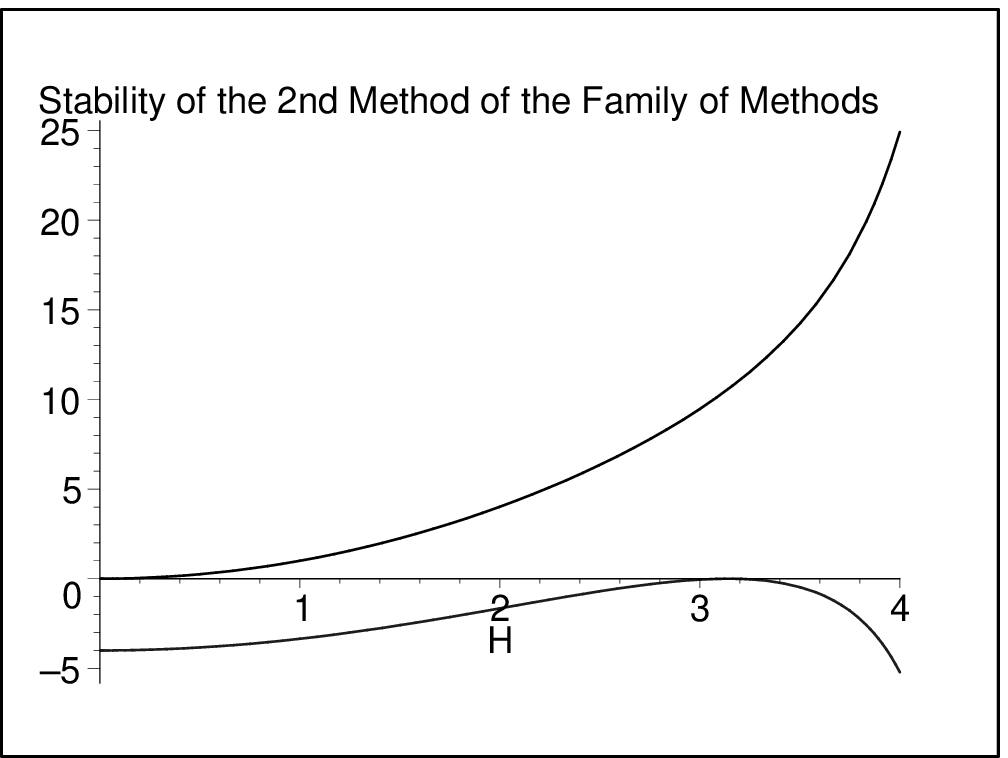}\\
\end{tabular}
\begin{tabular}{ccc}
\includegraphics[width=0.4\textwidth]{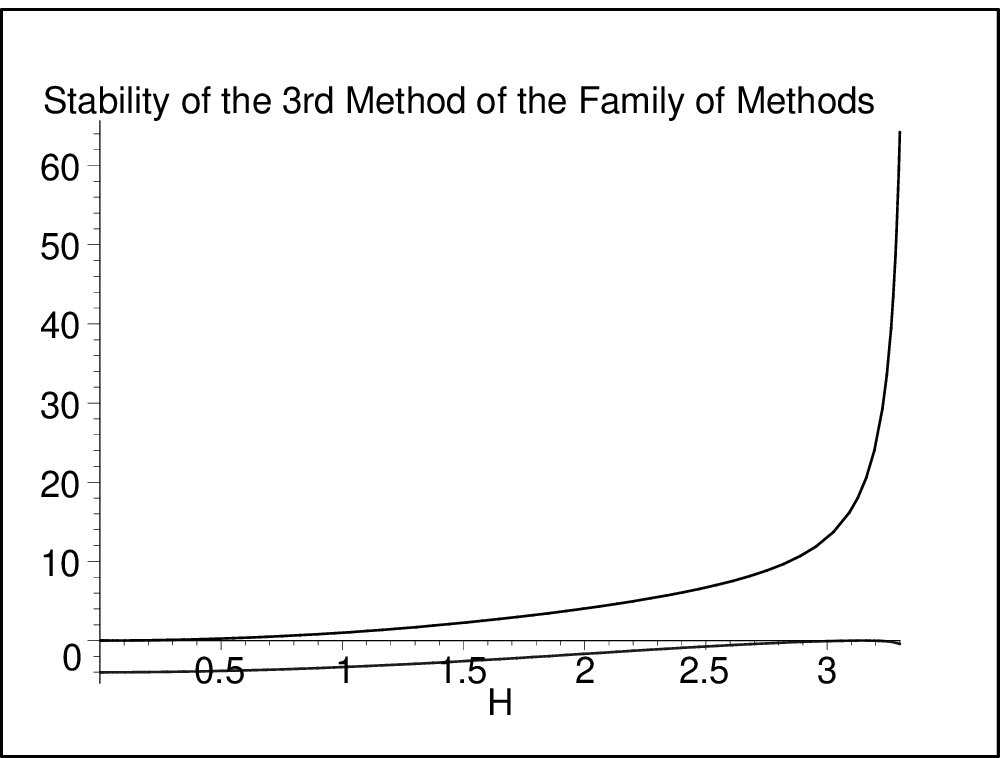}&\quad &
\includegraphics[width=0.4\textwidth]{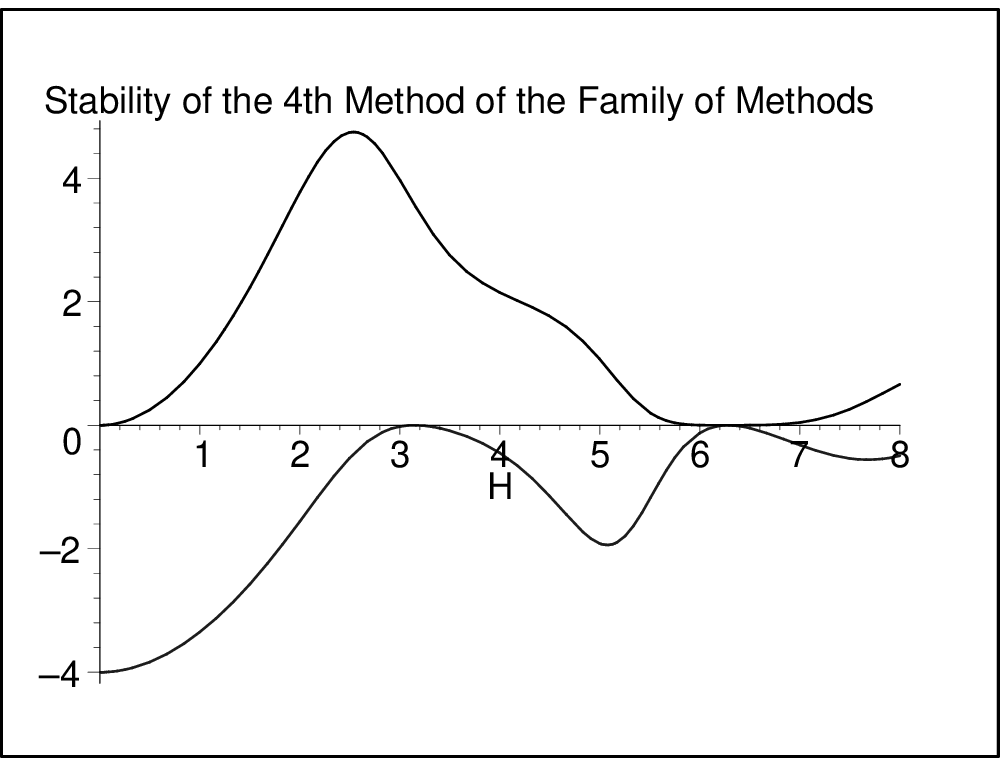}\\
\end{tabular}
\caption[]{\label{figs5}
Stability polynomials of the new developed methods in the case that
$H=s$} \normalsize
\end{center}
\end{figure}

\vspace*{0.7cm}

In the case that the frequency of the scalar test equation is equal
with the frequency of phase fitting, i.e. in the case that $H=s$, we
have the following figure for the stability polynomials of the new
developed methods. A method is P-stable if the $s-H$ plane is
not shadowed. From the above diagrams it is easy for one to
see that the interval of periodicity of all the new methods is equal
to: $\Bigl ( 0, \pi^2 \Bigr )$.

\begin{remark}
For the solution of the Schr\"odinger equation the frequency of
the exponential fitting is equal to the frequency of the scalar
test equation. So, it is necessary to observe the surroundings of
the first diagonal of the $w-H$ plane.
\end{remark}

\section{Numerical results - Conclusion}

In order to illustrate the efficiency of the new methods obtained in
paragraphs 3.1 - 3.4, we apply them to the radial time independent
Schr\"odinger equation.
\par
In order to apply the new methods to the radial Schr\"odinger
equation the value of parameter $v$ is needed. For every problem of
the one-dimensional Schr\"odinger equation given by (\ref{eqn:1})
the parameter $v$ is given by
\begin{equation}
v=\sqrt{|q(x)|}=\sqrt{|V(x)-E|} \label{freq}
\end{equation}
\noindent where $V(x)$ is the potential and $E$ is the energy.

\subsection{Woods-Saxon potential}
\par
We use the well known Woods-Saxon potential given by
\begin{equation}
V(x) = \frac{u_{0}}{1+z} - \frac{u_{0} z}{a \left ( 1 + z \right
)^{2} } \label{irkbound2}
\end{equation}
\noindent with $z = exp \left [ \left ( x - X_{0} \right )/a
\right ], \hbox{ } u_{0} = - 50, \hbox{ } a = 0.6$, and
$X_{0}=7.0$.
\par
The behavior of Woods-Saxon potential is shown in  the Figure 10.

\setcounter{figure}{9}
\begin{figure}[htb!]
\begin{center}
\begin{tabular}{ccc}
\includegraphics[width=0.4\textwidth]{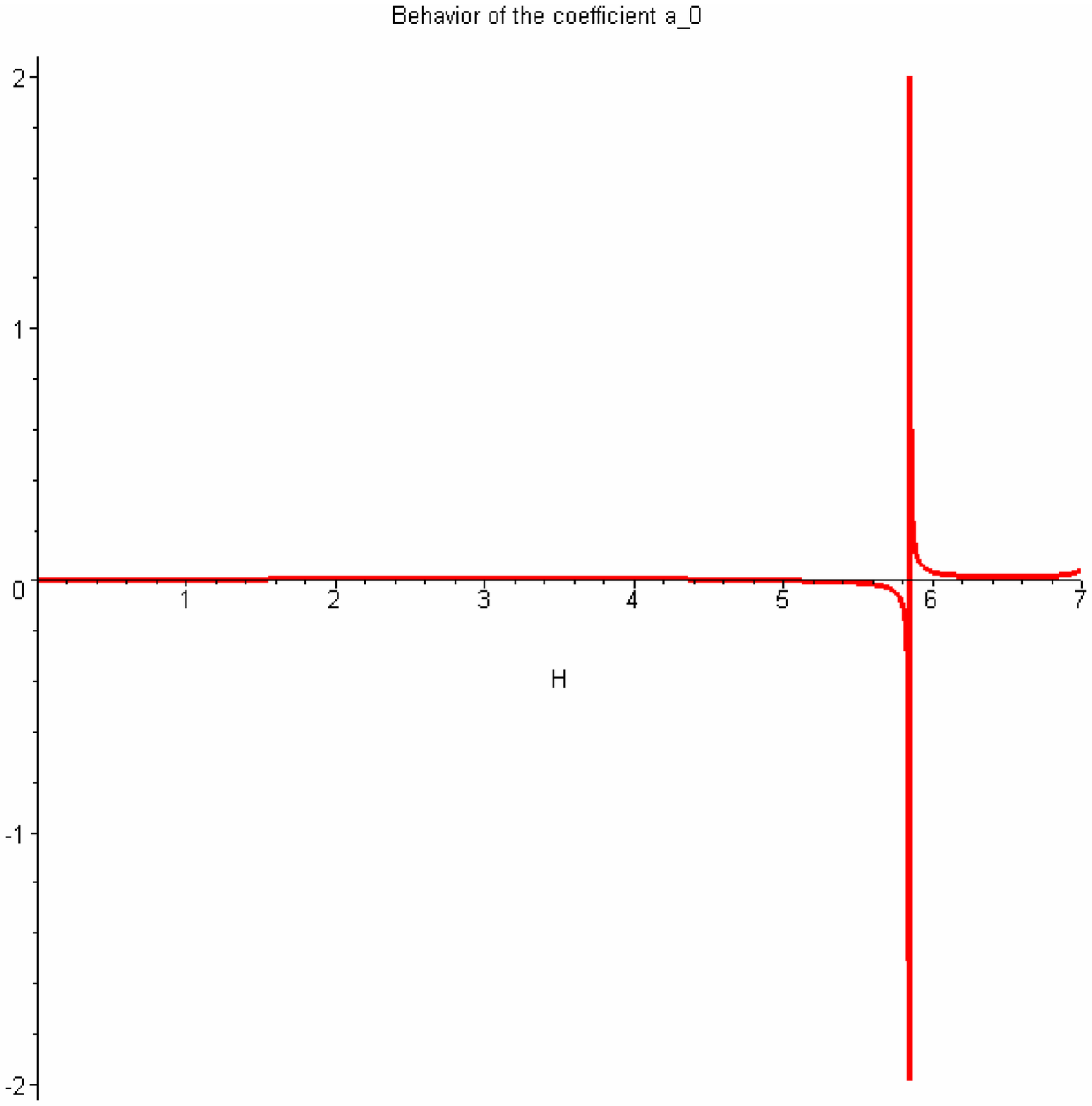}&\quad &
\includegraphics[width=0.4\textwidth]{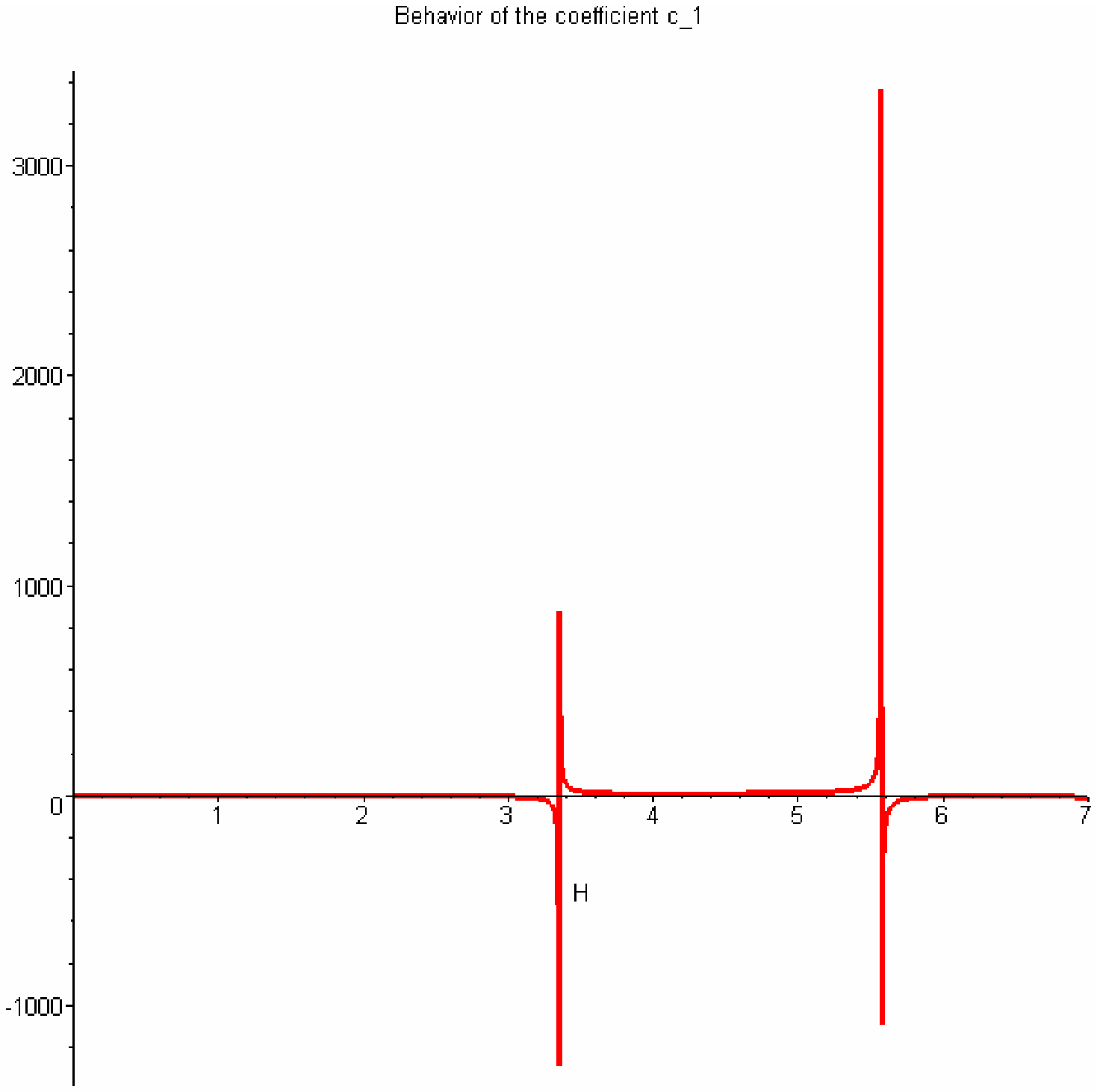}\\
\end{tabular}
\includegraphics[width=0.4\textwidth]{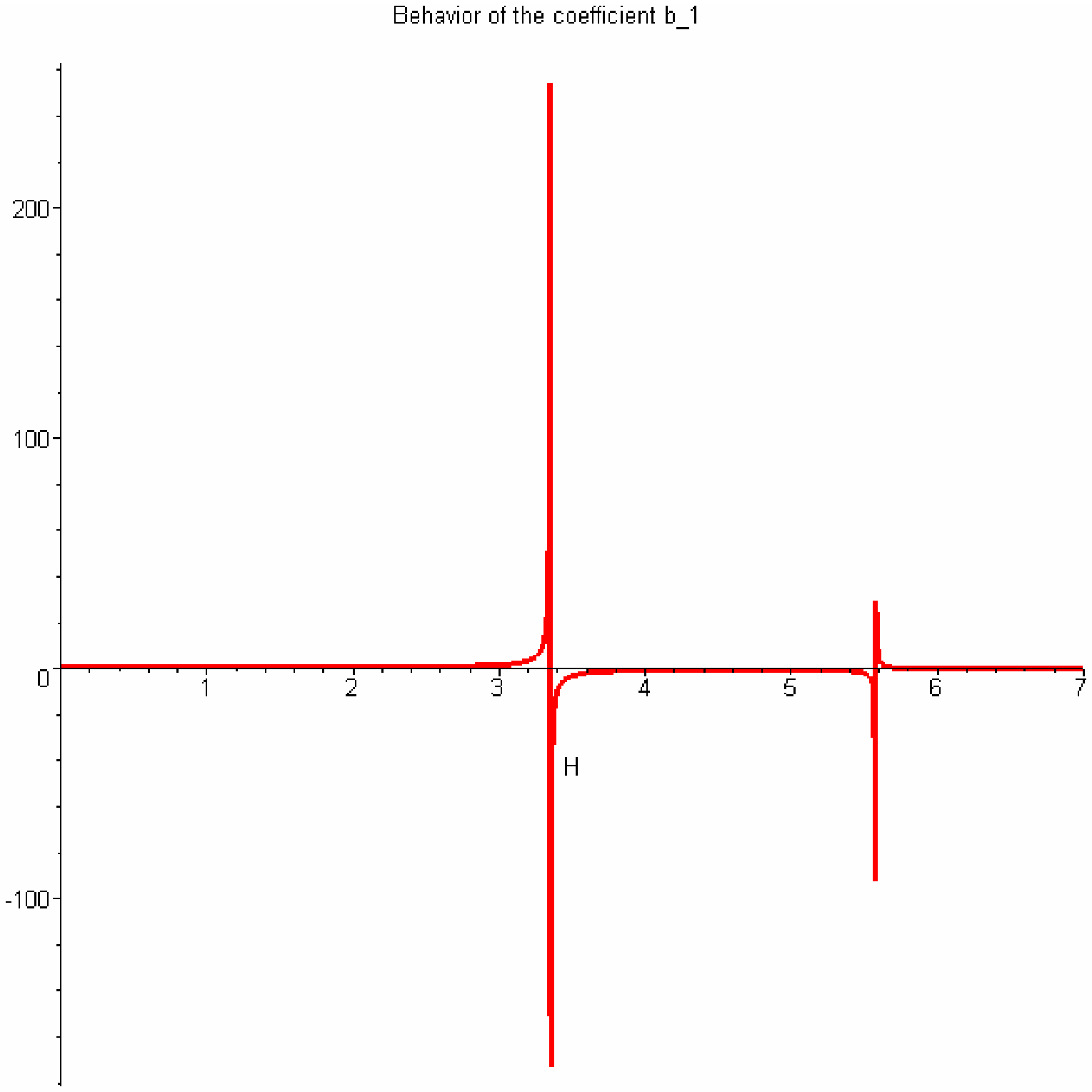}
\caption[]{\label{fig3}
Behavior of the coefficients of the new method given by (\ref{nm4b})
for several values of $H$.} \normalsize
\end{center}
\end{figure}

\vspace*{0.0cm}

\par
It is well known that for some potentials, such as the Woods-Saxon
potential, the definition of parameter $v$ is not given as a
function of $x$ but it is based on some critical points which have been
defined from the investigation of the appropriate potential (see for
details \cite{ix80}).

\par
For the purpose of obtaining our numerical results it is
appropriate to choose $v$ as follows (see for details
\cite{ix80}):

\begin{equation}
v = \left\{
\begin{array}{cl}
\sqrt{-50+E}, & \rm{for } \; x \in [0,6.5-2h],\\
\sqrt{-37.5+E}, & \rm{for } \; x = 6.5-\it{h}\\
\sqrt{-25+E} , & \rm{for } \; x = 6.5 \\
\sqrt{-12.5+E}, & \rm{for } \; x = 6.5+\it{h} \\
\sqrt{E}, & \rm{for } \; x \in [6.5+2h,15]
 \end{array}
 \right.
\end{equation}

\par
\subsection{Radial Schr\"odinger Equation - The Resonance Problem}
\par
Consider the numerical solution of the radial time independent
Schr\"odinger equation (\ref{eqn:1}) in the well-known case of the
Woods-Saxon potential (\ref{rkbound2}). In order to solve this
problem numerically we need to approximate the true (infinite)
interval of integration by a finite interval. For the purpose of
our numerical illustration we take the domain of integration as $x
\in [0,15]$.  We consider equation (\ref{eqn:1}) in a rather large
domain of energies, i.e. $E \in [1,1000]$.
\par
In the case of positive energies, $E=k^{2}$, the potential dies
away faster than the term $\frac{l(l+1)}{x^{2}}$ and the
Schr\"odinger equation effectively reduces to
\begin{equation}
y''(x) + \left ( k^{2} - \frac{l(l+1)}{x^{2}} \right ) y(x) = 0
\label{app1}
\end{equation}
\noindent for $x$ greater than some value $X$.
\par
The above equation has linearly independent solutions
$kxj_{l}(kx)$ and $kxn_{l}(kx)$  where $j_{l}(kx)$ and $n_{l}(kx)$
are the spherical Bessel and Neumann functions respectively. Thus
the solution of equation (\ref{eqn:1}) (when $x \rightarrow
\infty$ ) has the asymptotic form
\begin{eqnarray}
y(x) \simeq A kxj_{l}(kx) - B kxn_{l}(kx) \nonumber \\
\simeq AC \left [ sin \left ( kx - \frac{l \pi}{2} \right ) + tan
\delta_{l} cos \left ( kx - \frac{l \pi}{2} \right ) \right ]
\label{app2}
\end{eqnarray}
\noindent where $\delta_{l}$ is the phase shift, that is
calculated from the formula
\begin{equation}
tan \delta_{l} = \frac{y(x_{2}) S(x_{1}) - y(x_{1}) S(x_{2})}
{y(x_{1}) C(x_{1}) - y(x_{2}) C(x_{2})} \label{app3}
\end{equation}
for $x_{1}$ and $x_{2}$ distinct points in the asymptotic region
(we choose $x_{1}$ as the right hand end point of the interval of
integration and $x_{2} = x_{1} - h$) with $S(x) = kxj_{l}(kx)$ and
$C(x) = - kxn_{l}(kx)$. Since the problem is treated as an
initial-value problem, we need $y_{0}$ before starting a one-step
method. From the initial condition we obtain $y_{0}$. With these
starting values we evaluate at $x_{1}$  of the asymptotic region
the phase shift $\delta_{l}$.
\par
For positive energies we have the so-called resonance problem.
This problem consists either of finding the phase-shift
$\delta_{l}$ or finding those $E$, for $E \in [1,1000]$, at which
$\delta_{l} = \frac{\pi}{2}$. We actually solve the latter
problem, known as {\bf the resonance problem} when the positive
eigenenergies lie under the potential barrier.
\par
The boundary conditions for this problem are:
\begin{equation}
y(0) = 0, \hbox{ } y(x) = \cos \left ( \sqrt{E} x \right ) \mbox{
for large} \hbox{ }x. \label{app4}
\end{equation}



\setcounter{figure}{10}
\begin{figure}[htb!]
\begin{center}
\includegraphics[scale=0.55]{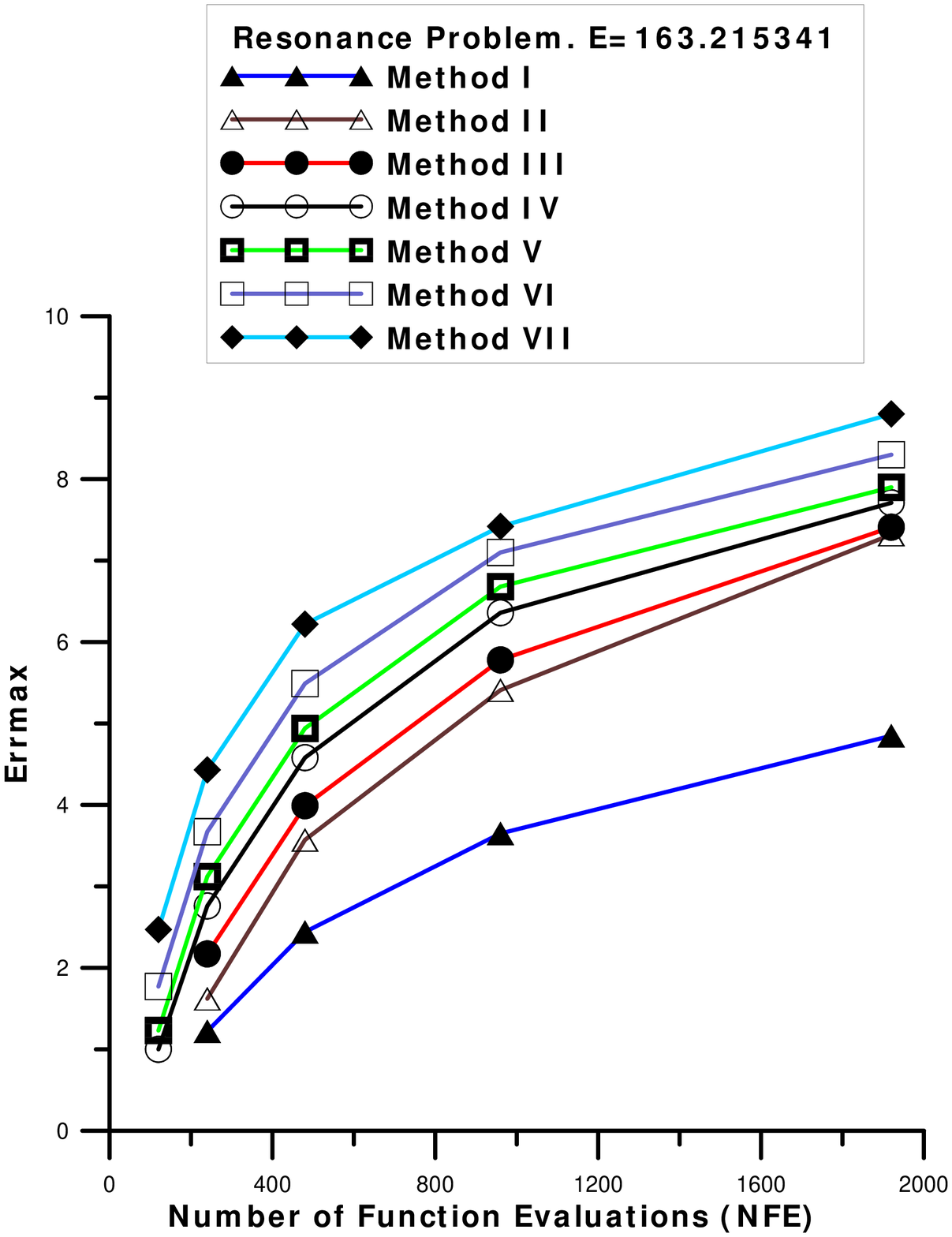}
\caption[]{\label{fig1aavve}
Error Errmax for several values of n for the eigenvalue $E_{1} =
163.215341$. The nonexistence of a value of Errmax indicates that
for this value of n, Errmax is positive} \normalsize
\end{center}
\end{figure}

\vspace*{0.0cm}

\setcounter{figure}{11}
\begin{figure}[htb!]
\begin{center}
\includegraphics[scale=0.55]{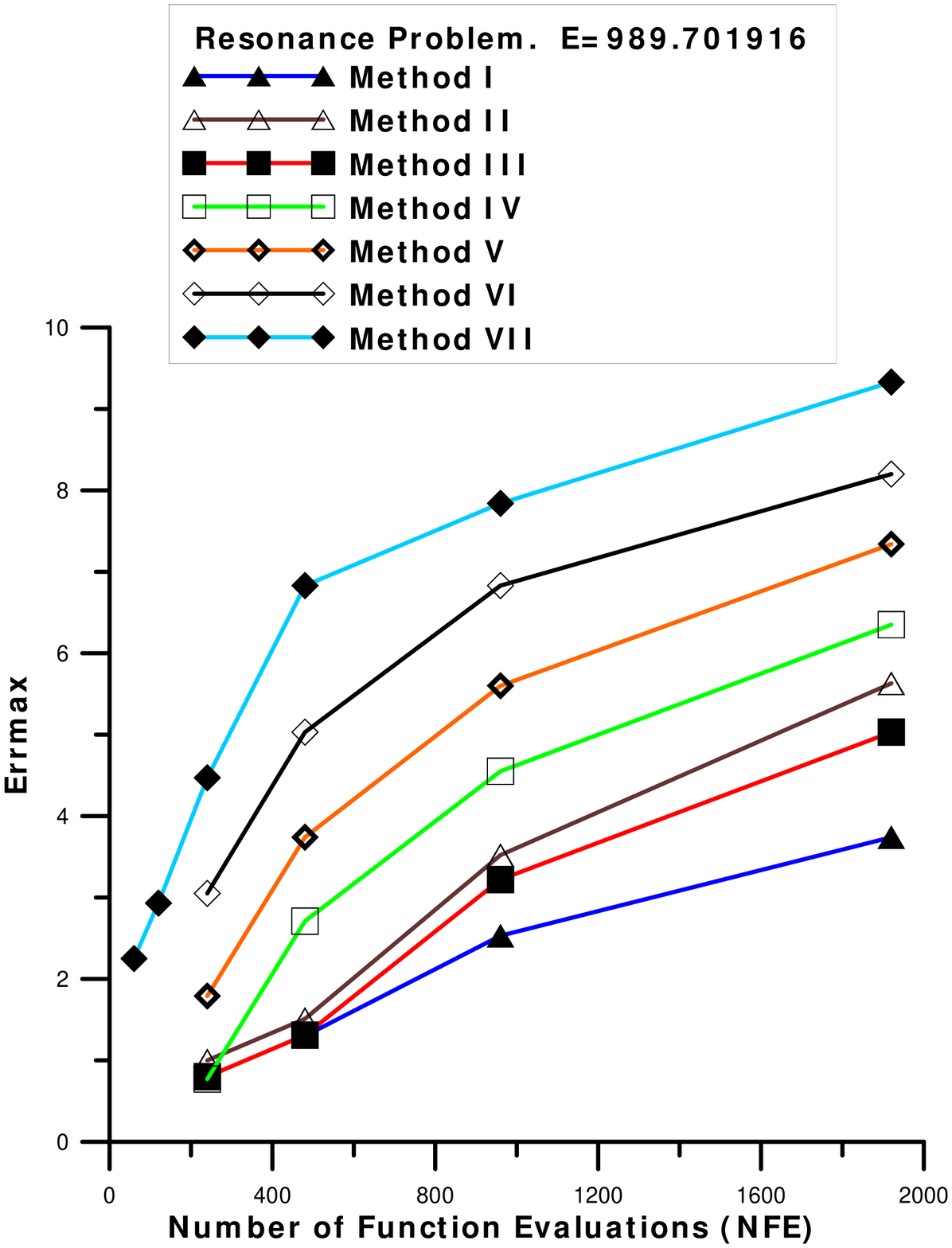}
\caption[]{\label{fig1aav}
Error Errmax for several values of n for the eigenvalue $E_{3} =
989.701916$. The nonexistence of a value of Errmax indicates that
for this value of n, Errmax is positive} \normalsize
\end{center}
\end{figure}

\vspace*{0.0cm}

\par
We compute the approximate positive eigenenergies of the
Woods-Saxon resonance problem using:

\begin{itemize}

\item The Numerov's method which is indicated as \textbf{Method I}

\item The Exponentially-fitted four-step method developed by Raptis \cite{ra83} which is indicated as
\textbf{Method II}

\item The Two-Step Numerov-type Method with minimum phase-lag produced by Chawla and Rao \cite{chawla84}  which is indicated as \textbf{Method III}

\item The new Two-Step Numerov-Type Method with phase-lag equal to zero obtained in paragraph 3.1 which is indicated as \textbf{Method IV}.

\item The new Two-Step Numerov-Type Method with phase-lag and its first derivative equal to zero obtained in paragraph 3.2 which is indicated as \textbf{Method V}.

\item The new Two-Step Numerov-Type Method with phase-lag and its first and second derivatives equal to zero obtained in paragraph 3.3 which is indicated as \textbf{Method VI}.

\item The new Two-Step Numerov-Type Method with phase-lag and its first, second and third derivatives equal to zero obtained in paragraph 3.4 which is indicated as \textbf{Method VII}.
\end{itemize}

The computed eigenenergies are compared with exact ones.  In Figure
11 we present the maximum absolute error $log_{10} \left ( Err
\right )$ where
\begin{equation}
Err=\left | E_{calculated} - E_{accurate} \right |
\end{equation}
\noindent of the eigenenergy $E_{1}$, for several values of
NFE = Number of Function Evaluations. In Figure 12 we present the maximum absolute error
$log_{10} \left ( Err \right )$ where
\begin{equation}
Err=\left | E_{calculated} - E_{accurate} \right |
\end{equation}
\noindent of the eigenenergy $E_{3}$, for several values of
NFE = Number of Function Evaluations.

\section{Conclusions}

In the present paper we have developed a family of methods of sixth
algebraic order for the numerical solution of the radial
Schr\"odinger equation.

More specifically we have developed:

\begin{enumerate}

\item A Two-Step Numerov-Type Method with phase-lag equal to zero

\item A Two-Step Numerov-Type Method with phase-lag and its first derivative equal to zero

\item A Two-Step Numerov-Type Method with phase-lag and its first and second derivatives equal to zero

\item A Two-Step Numerov-Type Method with phase-lag and its first, second and third derivatives equal to zero

\end{enumerate}

We have applied the new method to the resonance problem of the
radial Schr\"odinger equation.

Based on the results presented above we have the following
conclusions:

\begin{itemize}

\item The Exponentially-fitted four-step method developed by Raptis \cite{ra83} (denoted as Method II) is more efficient than the Numerov's Method (denoted Method I).

\item The Two-Step Numerov-type Method with minimum phase-lag produced by Chawla and Rao \cite{chawla84} (Method III) is more efficient than  the Exponentially-fitted four-step method developed by Raptis \cite{ra83} (Method II) for the energy
$163.215341$ and less efficient for the energy $989.701916$.

\item The new developed methods are much more efficient than the older ones.

\item The Two-Step Numerov-Type Method with phase-lag and its first derivative equal to zero (Method V) is more efficient than the New Two-Step Numerov-Type Method with phase-lag equal to zero
(Method IV)

\item The Two-Step Numerov-Type Method with phase-lag and its first and second derivatives equal to zero (Method VI) is more efficient than the Two-Step Numerov-Type Method with phase-lag and its first derivative equal to zero (Method V)

\item The Two-Step Numerov-Type Method with phase-lag and its first, second and third derivatives equal to zero (Method VII) is more efficient than the Two-Step Numerov-Type Method with phase-lag and its first and second derivatives equal to zero (Method VI)

\end{itemize}

\par
All computations were carried out on a IBM PC-AT  compatible 80486
using double precision arithmetic with 16 significant  digits
accuracy (IEEE standard).

\end{article}
\end{document}